\documentclass[10pt]{amsart}
\usepackage{amsthm, amssymb, amscd, mathptmx, pstricks}

\DeclareSymbolFont{bbold}{U}{dsrom}{m}{n}
\DeclareSymbolFontAlphabet{\mathbb}{bbold}

\renewcommand{\Bbb}{\mathbb}
\renewcommand{\frak}{\mathfrak}
\newcommand{\cal}{\mathcal}
\newcommand{\catqot}{/\hskip-3pt/}
\newcommand{\C}{{\Bbb C}}

\newcommand{\E}{{\cal E}}
\newcommand{\e}{\Bbb E}

\newcommand{\End}{\mathop{\rm End}}

\newcommand{\f}{\Bbb F}
\newcommand{\F}{{\cal F}}
\newcommand{\formall}{[\hskip-1.5pt[}
\newcommand{\formalr}{]\hskip-1.5pt]}
\newcommand{\G}{{\cal G}}
\newcommand{\PGL}{\mathop{\rm PGL}}
\newcommand{\GL}{\mathop{\rm GL}}

\newcommand{\Hom}{\mathop{\rm Hom}}

\renewcommand{\Im}{\mathop{\rm Im}}

\newcommand{\id}{\mathop{\rm id}}

\renewcommand{\L}{{\cal L}}

\newcommand{\M}{{\cal M}}

\newcommand{\n}{{\cal N}}
\newcommand{\Oh}{{\cal O}}

\newcommand{\Pe}{{\Bbb P}}

\newcommand{\Q}{{\Bbb Q}}
\newcommand{\R}{{\Bbb R}}

\newcommand{\SL}{\mathop{\rm SL}}
\newcommand{\Spec}{\mathop{\rm Spec}}

\newcommand{\trk}{\mathop{\rm trk}}

\newcommand{\Z}{{\Bbb Z}}

\newcommand{\la}{\lambda}
\newcommand{\lra}{\longrightarrow}

\newcommand{\lma}{\longmapsto}
\newcommand{\p}{\prime}
\newcommand{\q}{\quad}

\renewcommand{\phi}{\varphi}
\newcommand{\rk}{\mathop{\rm rk}}
\newcommand{\eps}{\varepsilon}

\renewcommand{\theta}{\vartheta}
\newcommand{\ul}{\underline}
\newcommand{\ol}{\overline}

\theoremstyle{plain}
\newtheorem{Thm}{Theorem}[subsection]
\newtheorem*{Thm*}{Main Theorem}
\newtheorem{Cor}[Thm]{Corollary}
\newtheorem{Prop}[Thm]{Proposition}
\newtheorem*{Prop*}{Proposition}
\newtheorem{Lem}[Thm]{Lemma}

\theoremstyle{definition}
\newtheorem{Def}[Thm]{Definition}

\theoremstyle{remark}
\newtheorem{Obs}[Thm]{Observation}
\newtheorem{Rem}[Thm]{Remark}
\newtheorem*{Rem*}{Remark}
\newtheorem{Ass}[Thm]{Assumption}
\title{The Hilbert compactification of the universal moduli space of semistable vector bundles
over smooth curves$^\dagger$\footnote{$^\dagger$to appear in J.\ Differential Geometry}}
\author{Alexander Schmitt}
\address{Universit\"at Duisburg-Essen, Fachbereich 6: Mathematik und Informatik, D-45117 Essen, Deutschland}
\email{alexander.schmitt@uni-essen.de}
\begin{document}


\maketitle
\markboth{Alexander Schmitt}{Hilbert Compactification}

\begin{abstract}
We construct the Hilbert compactification of the universal moduli space of semistable
vector bundles over smooth curves. The Hilbert compactification is the GIT quotient of
some open part of an appropriate Hilbert scheme of curves in a Gra\ss mannian.
It has all the properties asked for by Teixidor.
\end{abstract}

\section*{Introduction}
For every smooth curve $C$ and integers $\chi$ and $r>0$, one has
the projective moduli space ${\cal U}(C;\chi,r)$ of semistable
vector bundles $\E$ of rank $r$ and Euler characteristic
$\chi(\E)=\chi$. An automorphism $\sigma$ of $C$ acts on ${\cal
U}(C;\chi,r)$ via $[\E]\lma [\sigma^*\E]$. Let ${\frak M}_g$ be
the moduli space of smooth curves of genus $g$. It is possible to
construct a universal moduli space ${\frak U}(g;\chi,r)\lra
{\frak M}_g$, such that the fibre over $[C]$ is ${\cal
U}(C;\chi,r)/{\rm Aut}(C)$. This leads to the problem of
compactifying ${\frak U}(g;\chi,r)$ over $\ol{\frak M}_g$, the
moduli space of stable curves of genus $g$. There are two natural
approaches to this \cite{Teix}: First, given a stable curve $C$
of genus $g$, one can look at torsion free sheaves $\E$ of
uniform rank $r$ on $C$ with Euler characteristic $\chi(\E)=\chi$
which are semistable w.r.t.\ the canonical polarization. These
objects form again a projective moduli space ${\cal
U}(C;\chi,r)$. Pandharipande \cite{Pa} has constructed a
projective moduli space $\ol{\frak U}(g;\chi,r)\lra \ol{\frak
M}_g$, such that the fibre over a stable curve $[C]$ is ${\cal
U}(C;\chi,r)/{\rm Aut}(C)$. Second, for a stable curve $C$,
instead of looking at torsion free sheaves on $C$, one can look
at vector bundles on semistable models of $C$. This viewpoint has
advantages for certain degeneration arguments. As an approach to
the above problem, it has been formalized by Gieseker \cite{Gie}
and further studied by Gieseker and Morrison \cite{GM}, Nagaraj
and Seshadri \cite{Sesh}, Teixidor i Bigas \cite{Teix}, Kausz
(\cite{Ka1},\cite{Ka2}), and the author \cite{Sch}. It was also
used by Caporaso \cite{Cap} to solve the problem for $r=1$.
Without loss of generality, we may assume that, for every smooth
curve of genus $g$ and every semistable vector bundle $\E$ on $C$
of rank $r$ with $\chi(\E)=\chi$, $\E$ is globally generated,
$H^1(\E)=0$, and the evaluation map ${\rm ev}\colon
H^0(\E)\otimes{\Oh}_C\lra\E$ gives rise to a closed embedding
$C\hookrightarrow {\rm Gr}(H^0(\E), r)$ into the Gra\ss mannian
of $r$-dimensional quotients of $H^0(\E)$. Thus, we fix a vector
space $V^\chi$ of dimension $\chi$, define ${\frak G}:= {\rm
Gr}(V^\chi,r)$, and look at ${\frak H}(g;\chi,r)$, the closure of
the Hilbert scheme of smooth curves in ${\frak G}$ with Hilbert
polynomial $P(m)=d\cdot m+(1-g)$ in the whole Hilbert scheme.
Here, $\chi=d+r(1-g)$ and ${\frak G}$ is polarized by ${{\Oh}}_{\frak
G}(1)$, the determinant of the universal quotient bundle. Note
that we have a natural action of ${\rm SL}(V^\chi)$ on ${\frak
H}(g;\chi,r)$. Our candidate for the Hilbert compactification is,
therefore, ${\frak H\frak C}(g;\chi,r):={\frak
H}(g;\chi,r)\catqot \SL(V^\chi)$. Before we can form the
GIT-quotient, we have, however, to find appropriate linearized
ample line bundles on ${\frak H}(g;\chi,r)$. First, there are
some obvious ones. For this, let ${\frak C}\subset {\frak
G}\times {\frak H}(g;\chi,r)$ be the universal curve. For every
natural number $m$, we have, on ${\frak G}\times {\frak
H}(g;\chi,r)$, the exact sequence
$$
\begin{CD}
0 @>>> {\cal I}_{\frak C}\otimes \pi_{\frak G}^*{\Oh}_{\frak G}(m) @>>>
\pi_{\frak G}^*{\Oh}_{\frak G}(m) @>>> \bigl(\pi_{\frak G}^*{\Oh}_{\frak G}(m)\bigr)_{|\frak C}
@>>> 0.
\end{CD}
$$
For large $m$, $\pi_{{\frak H}(g;\chi,r)*}$ of this sequence leads to a surjective
homomorphism of vector bundles
$$
\begin{CD}
\Psi_{\frak H(g;\chi,r)}\colon S^m V^\chi \otimes{\Oh}_{\frak H(g;\chi,r)}
\lra {\frak E}_m:=\pi_{{\frak H}(g;\chi,r)*}
\Bigl(
\bigl(\pi_{\frak G}^*{\Oh}_{\frak G}(m)\bigr)_{|\frak C}\Bigr).
\end{CD}
$$
The rank of ${\frak E}_m$ is $P(m)$, and $\wedge^{P(m)} \Psi_{\frak H(g;\chi,r)}$
yields a closed immersion
$$
{\frak H}(g;\chi,r)\hookrightarrow {\Pe}\Bigl(\bigwedge^{P(m)} S^m V^\chi\Bigr)
$$
which
is equivariant w.r.t.\ the $\SL(V^\chi)$-actions. Let ${\frak L}_m:=
{\frak L}_m(g;\chi,r):=\det({\frak E}_m)$ be the pullback of
${\Oh}(1)$. Then, we can define ${\frak H\frak C}_m(g;\chi,r):={\frak H}\catqot_{\frak L_m}
\SL(V^\chi)$.
This, space will be, however, only useful, if we can settle the following properties:
\begin{enumerate}
\item If $[C]\in {\frak H}(g;\chi,r)$ is a smooth curve and $q_C\colon
V^\chi\otimes{\Oh}_C\lra \E$ is
the pullback of the universal quotient, then $[C]$ is (semi)stable w.r.t.\ the linearization
in ${\frak L}_m$, if and only if $\E$ is a (semi)stable vector bundle.
\item If $[C]\in {\frak H}(g;\chi,r)$ is semistable w.r.t.\ the linearization
in ${\frak L}_m$ and $q_C\colon V^\chi\otimes{\Oh}_C\lra \E$ is the pullback of the universal quotient,
then $C$ is a semistable curve and $V^\chi\lra H^0(\E)$ is an isomorphism.
\end{enumerate}
Point (1) and (2) have been settled in the rank two case by
Gieseker and Morrison \cite{GM}, and (1) in general by the author
\cite{Sch}. Unfortunately, nothing is known about (2) in general,
and, even if it were true, the computations of the correct notion
of semistability would still be extremely difficult (cf.\
\cite{Teix}). The way out is to adapt a strategy due to Nagaraj
and Seshadri \cite{Sesh}. In our setting, it is described as
follows: For every stable curve $C$, let ${\cal Q}(C;\chi,r)$ be
the quot scheme parameterizing quotients $V^\chi\otimes
{\Oh}_C\lra\E$ where $\E$ is a coherent sheaf of uniform rank $r$
with Euler characteristic $\chi(\E)=\chi$. From Pandharipande's
construction, we get a universal quot scheme $\ol{\frak
Q}(g;\chi,r)\lra \ol{\frak M}_g$, such that the fibre over $[C]$
is just ${\cal Q}(C;\chi,r)/{\rm Aut}(C)$, and a natural
$\SL(V^\chi)$-linearized ample line bundle ${\frak N}$ on
$\ol{\frak Q}(g;\chi,r)$. Next, let ${\frak H}^0(g;\chi,r)\subset
{\frak H}(g;\chi,r)$ be the open part corresponding to semistable
curves with the following property: If $\pi\colon C\lra C^\p$ is
the projection onto the stable model and if $q_C\colon
V^\chi\otimes{\Oh}_C\lra \E$ is the pullback of the universal
quotient on ${\frak G}$, then $\pi_*(q_C)\colon
V^\chi\otimes{\Oh}_{C^\p}\lra \pi_*(\E)$ is surjective and
$\pi_*(\E)$ is a torsion free sheaf which is semistable w.r.t.\
the canonical polarization. There is a natural map ${\frak
H}^0(g;\chi,r)\lra \ol{\frak Q}\ (g;\chi,r)$, landing in the
${\frak N}$-semistable locus. Let $\widehat{\frak
H}(g;\chi,r)\subset {\frak H}(g;\chi,r)\times \ol{\frak
Q}(g;\chi,r)$ be the closure of the graph of the above morphism.
This is an $\SL(V^\chi)$-invariant subscheme. For $a\gg 0$, the
semistable points w.r.t.\ the linearization in
$$
{\frak L}_{a,m}:={\frak L}_{a,m}(g;\chi,r)
:=\bigl(\pi_{{\frak H}(g;\chi,r)}^*{\frak L}_m(g;\chi,r)\otimes \pi_{\ol{\frak Q}(g;\chi,r)}^*
{\frak N}^{\otimes a}\bigr)_{|\widehat{\frak H}(g;\chi,r)}
$$
will lie in the set of the preimages of the points in $\ol{\frak Q}(g;\chi,r)$
which are semistable w.r.t.\ the linearization in ${\frak N}$.
Note that, for every $l\ge 0$, we can perform the same constructions w.r.t.\
$g$, $r$, and $\chi_l:=\chi+l\cdot r\cdot (2g-2)$. The result of this note is
\begin{Thm*}
There exist an $l_0$ and for every $l\ge l_0$ an $m(l)$, such that
for all $l\ge l_0$ and $m\ge m(l)$ the following
properties hold true:
\par
{\rm i)} All points in
$\widehat{\frak H}(g;\chi_l,r)$ which are semistable w.r.t.\ the linearization in
${\frak L}_{a,m}(g;\chi_l,r)$ for all $a\gg 0$ lie in the graph
of the morphism ${\frak H}^0(g;\chi_l,r)\lra \ol{\frak Q}(g;\chi_l,r)$.
\par
{\rm ii)} Let $[C]$ be contained in the graph of the
morphism ${\frak H}^0(g;\chi_l,r)\lra \ol{\frak Q}(g;\chi_l,r)$, and let
$q_C\colon V^{\chi_l}\otimes{\Oh}_C\lra \E$ be the pullback of the universal quotient.
Then, $[C]$
is (semi)stable w.r.t.\ the linearization in
${\frak L}_{a,m}(g;\chi_l,r)$ for all $a\gg 0$, if and only if
$(C,\E\otimes\omega_C^{\otimes -l})$ is
H-(semi)stable in the sense of Definition~{\rm\ref{H-stab}}.
\end{Thm*}
\begin{Rem*}
i) If $\E$ is a vector bundle of rank $r$ on the semistable curve
$C$ and if $\pi\colon C\lra C^\p$ is the map onto the stable
model, then the condition that $\pi_*(\E)$ be torsion free is a
precise condition on the restriction of $\E$ to any chain $R$ of
rational curves attached at only two points, say $p_1$ and $p_2$.
Namely, there must not exist a non-zero section of $\E_{|R}$ vanishing in
both $p_1$ and $p_2$. This implies, for example, that $R$ has at
most $r$ components and that $\E_{|R}$ is \it strictly standard\rm,
i.e., for any component $R_i\cong{\Pe}_1$ of $R$, $\E_{|R_i}\cong
{\Oh}_{R_i}^{\oplus \eps_i}\oplus {\Oh}_{R_i}(1)^{\oplus (r-\eps_i)}$
with $0\le \eps_i<r$, $i=1,...,s$. For the detailed discussion,
we refer the reader to \cite{Sesh}.
\par
ii) The condition to be semistable w.r.t.\ the linearization in
${\frak L}_{a,m}(g;\chi_l,r)$ \sl for all $a\gg 0$ \rm is
explained as follows: Let $G$ be a reductive algebraic group,
$\rho_i\colon G\lra \GL(W_i)$, $i=1,2$, two finite dimensional
representations, and ${\widehat{\frak H}}\subset
{\Pe}(W_1)\times{\Pe}(W_2)$ a $G$-invariant closed subscheme. Denote by
${\frak L}_{a,m}$ the restriction of
${\Oh}_{{\Pe}(W_1)\times{\Pe}(W_2)}(m,a)$ to ${\widehat{\frak H}}$ and by
${\widehat{\frak H}}_{a,m}^{\rm (s)s}$ the set of points which are
(semi)stable w.r.t.\ the linearization in ${\frak L}_{a,m}$.
Then, there is an $\eps_\infty$, such that
$$
\pi_2^{-1}\bigl({\Pe}(W_2)^{\rm s}\bigr)\q \subset\q{\widehat{\frak
H}}^{\rm (s)s}_{a^\p,m^\p}\q=\q {\widehat{\frak H}}^{\rm (s)s}_{a,m}
\q\subset\q \pi_2^{-1}\bigl({\Pe}(W_2)^{\rm ss}\bigr)
$$
whenever both $m^\p/a^\p$ and $m/a$ are smaller than
$\eps_\infty$. This follows from the corresponding assertion for
$\C^*$-actions and the master space construction of Thaddeus
(\cite{Th}, \cite{OST}).
\par
iii) The intrinsically defined concept of H-(semi)stability has
the following properties:
\begin{itemize}
\item For smooth curves, it agrees with Mumford-(semi)stability.
\item If $(C,\E\otimes\omega_C^{\otimes -l})$ is H-semistable and
$\pi\colon C\lra C^\p$ is the morphism to the stable model, then
$\pi_*(\E\otimes\omega_C^{\otimes -l})$ is semistable w.r.t.\ the
canonical polarization.
\item A pair $(C,\E\otimes\omega_C^{\otimes -l})$ with $C$ a semistable curve and
$\pi_*(\E\otimes\omega_C^{\otimes -l})$ a stable sheaf,
$\pi\colon C\lra C^\p$ being the morphism to the stable model, is
H-stable.
\end{itemize}
\end{Rem*}
Therefore, ${\frak H\frak C}(g;\chi,r):=\widehat{\frak H}(g;\chi_l,r)\catqot_{{\frak L}_{a,m}
(g;\chi_l,r)}\SL(V^{\chi_l})$
with $l\ge l_0$, $m\ge m(l)$, and $a\gg 0$ is a well-defined moduli space which compactifies
${\frak U}(g;\chi,r)$ over $\ol{\frak M}_g$ and, furthermore, maps to $\ol{\frak U}(g;\chi,r)$.
We remark that the authors of \cite{Sesh} were well aware of the
fact that their approach might be used in this generality. The
formulation of an intrinsic semistability concept for vector
bundles on semistable curves, however, seems to be new.
\par
The final section of this paper is devoted to the study of the
geometry of the Hilbert compactification and its map to the
moduli space of stable curves.
\subsection*{Acknowledgment}
The author wishes to thank the referee for some valuable
suggestions and Prof.\ Newstead for pointing out reference
\cite{Aut} to him.
\par
During the preparation of the article, the author benefited from support by the DFG through
the priority program
``Globale Methoden in der komplexen Geometrie --- Global Methods in Complex Geometry''.
\subsection*{Notation}
Let $C$ be a semistable curve. Then, ${\Oh}_C(1)$ will stand for the
canonical sheaf $\omega_C$, although this line bundle will be
ample if and only if $C$ is stable. Likewise, if $\E$ is a
coherent sheaf on $C$, we write $P(\E)$ for the polynomial $l\lma
\chi(\E(l))$. A \it scheme \rm will be a scheme of finite type
over the field of complex numbers.
\section{Preliminaries}
\subsection{Suitable linearizations}
Let $G$ be a reductive group, and $\rho_i\colon G\lra \GL(W_i)$, $i=1,2$, two representations
of $G$. This yields an action of $G$ on ${\Pe}(W_1)\times{\Pe}(W_2)$. Assume ${\frak X}$
is a $G$-invariant subscheme, and let $\pi\colon {\frak X}\lra {\Pe}(W_2)$ be the induced morphism.
Finally, let ${\frak L}_{n_1,n_2}$ be the $G$-linearized ample line bundle
${\Oh}_{{\Pe}(W_1)\times{\Pe}(W_2)}(n_1,n_2)_{|{\frak X}}$.
Define ${\frak X}^{\rm (s)s}_{n_1,n_2}$ as the set of points in ${\frak X}$
which are (semi)stable w.r.t.\ the linearization in
${\frak L}_{n_1,n_2}$. Likewise, ${\Pe}(W_2)^{\rm (s)s}$ is defined. The following
is well known and easy to prove.
\begin{Prop}
\label{suitable}
If $n_2/n_1$ is large enough, then
$$
\pi^{-1}\bigl({\Pe}(W_2)^{\rm s}\bigr)\q\subset\q{\frak X}_{n_1,n_2}^{\rm s}
\q\subset\q {\frak X}_{n_1,n_2}^{\rm ss}\q\subset\q \pi^{-1}\bigl({\Pe}(W_2)^{\rm ss}\bigr).
$$
\end{Prop}
\begin{Rem}
\label{Absi}
Let $n_2/n_1$ be so large that the conclusion of \ref{suitable} holds.
Then, a point $x\in {\frak X}$ will be (semi)stable w.r.t.\ the linearization in
${\frak L}_{n_1,n_2}$, if and only if $\pi(x)$ is semistable and for every
one parameter subgroup $\la\colon \C^*\lra G$ with $\mu_{{\Oh}_{{\Pe}(W_2)}(1)}(\la,\pi(x))=0$,
one has $\mu_{{\Oh}_{{\Pe}(W_1)}(1)}(\la,\pi^\p(x))\ (\ge)\ 0$, $\pi^\p\colon {\frak X}\lra{\Pe}(W_1)$
being the induced morphism.
\end{Rem}
\subsection{Pandharipande's moduli space and the universal quot scheme}
\label{Pandhi}
Let $C$ be a stable curve with irreducible components $C_1,...,C_c$, and $\F$ a coherent
sheaf on $C$. The tuple $\ul{r}(\F):=(\rk \F_{|C_1},...,\rk\F_{|C_c})$ is called
the \it multirank of $\F$\rm. We say that $\F$ has \it uniform rank $r$ on $C$\rm,
if $\ul{r}(\F)=(r,...,r)$. Finally, set $e_\iota:=\deg \omega_{C|C_\iota}$, $\iota=1,...,c$.
Then, the \it total rank of $\F$ \rm is the quantity $\trk\F:=\sum_{\iota=1}^c e_\iota
\cdot\rk \F_{|C_\iota}$.
Now, a coherent ${\Oh}_C$-module $\E$ is called \it (semi)stable\rm,
if it is torsion free and for every subsheaf $0\subsetneq\F\subsetneq\E$, the inequality
$$
\frac{\chi(\F)}{\trk\F}\q (\le)\q \frac{\chi(\E)}{\trk\E}
$$
is satisfied. One also introduces the notions of \it S-equivalence \rm and \it polystability\rm.
\rm
Pandharipande studies the functor $\ul{\rm U}(g;\chi,r)$ which associates to every
scheme $S$ the set of equivalence classes of pairs $(\cal C_S,\E_{S})$,
consisting of a flat family $\pi\colon \cal C_S\lra S$ of stable curves and an $S$-flat
coherent
sheaf $\E_S$ on ${\cal C}_S$ such that $\E_{S|\pi^{-1}(s)}$ is a semistable
sheaf of uniform rank $r$ and Euler characteristic
$\chi$
on $\pi^{-1}(s)$ for all closed points $s\in S$. Here,
$(\cal C_S,\E_{S})$ and $(\cal C^\p_S,\E^\p_{S})$ are considered \it equivalent\rm,
if there are an $S$-isomorphism $\psi\colon {\cal C}_S\lra {\cal C}_S^\p$ and
a line bundle $\L_S$ on $S$ such that $\E_{S}\cong \psi^*\E^\p_S\otimes\pi^*\L_S$.
In \cite{Pa}, a coarse moduli space $\ol{\frak U}(g;\chi,r)$ for the functor
$\ul{\rm U}(g;\chi,r)$ is constructed.
\par
We now briefly review the construction, because we need some of the details.
If $C$ is a stable curve, then ${\Oh}_C(10)$ defines a closed embedding
$C\hookrightarrow {\Pe}_N$, $N:=10(2g-2)-g$. Let ${\frak H}_g$ be the Hilbert scheme
of curves in ${\Pe}_N$ with the respective Hilbert polynomial. There is a natural left
action of $\SL(N+1)$ on ${\frak H}_g$ together with a linearization in an ample
line bundle ${\frak L}_{\frak H_g}$. As Gieseker \cite{GieKurv} has shown, the
GIT quotient ${\frak H}_g\catqot_{\frak L_{\frak H_g}}\SL(N+1)$ yields $\ol{\frak M}_g$.
\par
There is a constant $l_0$, such that for every $l\ge l_0$, every stable curve
$C$ of genus $g$, and every (semi)stable sheaf $\E$ of uniform rank $r$
and Euler characteristic $\chi$ on $C$, one has
\begin{itemize}
\item $\E(10\cdot l)$ is globally generated.
\item $H^1(\E(10\cdot l))=0$.
\item After identification of $H^0(\E(10\cdot l))$ with a
      previously fixed vector space $V^{\chi_{10\cdot l}}$ of dimension $\chi_{10\cdot l}
      = \chi+10\cdot l\cdot r\cdot (2g-2)$, the point ${\rm ev}\colon V^{\chi_{10\cdot l}}
      \otimes {\Oh}_C\lra\E(10\cdot l)$ in the quot scheme ${\cal Q}(C;\chi_{10\cdot l},r)$
      is (semi)stable.
\end{itemize}
Let $\cal C_g\hookrightarrow {\frak H}_g\times{\Pe}_N$ be the universal curve. We then have
a relative quot scheme $\rho\colon Q=Q(g;\chi_{10\cdot l},r)\lra {\frak H}_g$, such that
the fibre over $C\in {\frak H}_g$  is the quot scheme ${\cal Q}(C;\chi_{10\cdot l},r)$
of quotients $\E$ of $V^{\chi_{10\cdot l}}\otimes {\Oh}_C$ with $\chi(\E\otimes_{{\Oh}_C}{\Oh}_{{\Pe}_N}(l))
=\chi+10\cdot l\cdot r\cdot (2g-2)$ for all $l$.
The natural action of $\SL(V^{\chi_{10\cdot l}})\times\SL(N+1)$ on $Q$
is linearized in a suitable $\rho$-ample line bundle ${\frak L}_Q$. For any $a>0$,
define ${\frak L}_a:={\frak L}_Q\otimes\rho^*{\frak L}_{\frak H_g}^{\otimes a}$.
For large $a$, we have, by Proposition~\ref{suitable},
$Q^{\rm ss}_{\frak L_a}\subset \rho^{-1}({\frak H}_{g, {\frak L}_{\frak H_g}}^{\rm ss})$.
Then, $\ol{\frak U}(g;\chi,r):=Q\catqot_{\frak L_a}
\bigl(\SL(V^{\chi_{10\cdot l}})\times\SL(N+1)\bigr)$.
\par
Now, we can form this GIT quotient in two steps \cite{OST}: First, we divide
by the $\SL(N+1)$-action and then by the $\SL(V^{\chi_{10\cdot l}})$-action.
As ${\frak H}_{g, {\frak L}_{\frak H_g}}^{\rm s}=
{\frak H}_{g, {\frak L}_{\frak H_g}}^{\rm ss}$ (more precisely, this holds
on an appropriate closed subscheme, see \cite{GieKurv}, Proposition~2.0.0),
\ref{suitable} shows that a point $(C,q_C)\in Q$
is $\SL(N+1)$-semistable w.r.t.\ the linearization in ${\frak L}_{\frak H_g}$,
if and only if $C$ is stable, i.e.,
there is no condition on the quotient $q_C$. Set $\ol{\frak Q}(g;\chi_{10\cdot l},r)
:=Q\catqot_{\frak L_a}\SL(N+1)$. Let $\ul{\rm Q}(g;\chi_{10\cdot l},r)$ be the functor
which assigns to a scheme $S$ the set of equivalence classes of pairs
$(\cal C_S, q_S\colon V^{\chi_{10\cdot l}}\otimes{\Oh}_{\cal C_s}\lra \E_S)$ consisting
of a flat family $\pi\colon {\cal C}_S\lra S$ of stable curves and a quotient
$q_S$ onto an $S$-flat sheaf $\E_S$, such that $\chi(\E_{S|\pi^{-1}(s)})
=\chi_{10\cdot l}$ for every closed point $s\in S$.
Two families $(\cal C_S,q_S)$ and $(\cal C_S^\p,q_S^\p)$ will be considered \it equivalent\rm,
if there are an $S$-isomorphism $\psi\colon \cal C_S\lra \cal C_S^\p$ and
an isomorphism $\phi_S\colon \E_S\lra \psi^*\E_S^\p$ with
$\psi^*q_S^\p=\phi_S\circ q_S$.
The space $\ol{\frak Q}(g;\chi_{10\cdot l},r)$ obviously is the coarse moduli scheme
for the functor $\ul{\rm Q}(g;\chi_{10\cdot l},r)$. In particular, the fibre over
$[C]\in\ol{\frak M}_g$ identifies with ${\cal Q}(C;\chi_{10\cdot l},r)/{\rm Aut}(C)$.
As explained in the first section of \cite{OST}, there is an induced
$\SL(V^{\chi_{10\cdot l}})$-action on $\ol{\frak Q}(g;\chi_{10\cdot l},r)$
and some multiple of ${\frak L}_a$
descends to an $\SL(V^{\chi_{10\cdot l}})$-linearized ample line bundle ${\frak N}$.
Moreover, the points in $\ol{\frak Q}(g;\chi_{10\cdot l},r)$ which are
$\SL(V^{\chi_{10\cdot l}})$-(semi)stable for the given linearization are just the images
of the $(\SL(V^{\chi_{10\cdot l}})\times\SL(N+1))$-(semi)stable points.
Let us note the following elementary fact.
\begin{Lem}
\label{equal}
Let $C$ be a stable curve and $q_C\colon V^{\chi_{10\cdot l}}\otimes{\Oh}_C\lra \E$
a quotient with $\chi(\E)=\chi_{10\cdot l}$.
Let $h\in {\frak H}_g$ be a point such that the fibre of the universal curve
over $h$ is isomorphic to $C$.
Then, for any one parameter subgroup $\la\colon \C^*\lra \SL(V^{\chi_{10\cdot l}})$,
$$
\mu_{{\frak L_a}}\bigl(\la, (h,q_C)\bigr)\ (\ge)\ 0
\q\Longleftrightarrow\q
\mu_{{\frak N}}\bigl(\la, [C,q_C]\bigr)\ (\ge)\ 0.
$$
\end{Lem}
\begin{proof}
The quotient map $Q^{\SL(N+1)-\rm ss}_{\frak L_a}\lra
\allowbreak
\ol{\frak Q}(g;\chi_{10\cdot l},r)$
yields the $\SL(V^{\chi_{10\cdot l}})$-equivariant and finite morphism
$Q_h:=\rho^{-1}(h)\lra\allowbreak {\cal Q}(C;\chi_{10\cdot l},r)/{\rm Aut}(C)$. This immediately implies the assertion.
\end{proof}
\subsection{Auxiliary results from Nagaraj and Seshadri}
We recall some of the results of the paper \cite{Sesh} which we will use. Additional information
may be found there. We also point out that the above paper works
with semistable curves the stable model of which is an irreducible curve with exactly one
node. As one easily check, this assumption is not essential.
\begin{Prop}
\label{NSproperties}
{\rm i)} Let $R$ be a chain of projective lines, and $\F$ a globally generated torsion
free sheaf on $R$. Then, for any component $R_\iota\cong{\Pe}_1$, the restriction map
$H^0(R,\F)\lra H^0(R_\iota,\F_{|R_\iota})$ is surjective. Moreover, $H^1(R,\F)=0$.
\par
{\rm ii)} Let $\pi\colon C^\p\lra C$ be a morphism between semistable curves which
contracts only some chains of projective lines. Suppose $\E$ is a vector bundle
on $C^\p$ the restriction of which to every projective line contracted by $\pi$
has non-negative degree. In that situation
\begin{enumerate}
\item $\pi_*{\Oh}_{C^\p}={\Oh}_C$.
\item $R^i\pi_*(\E)=0$ for $i>0$. In particular, $H^j(C^\p,\E)=H^j(C,\pi_*(\E))$ for
all $j$.
\item Let $R^j$ be a chain of projective lines which is contracted by $\pi$
and attached at the points
$p_1$ and $p_2$. Let $\widetilde{C}$ be the closure of $C^\p\setminus R^j$.
If $H^1(\widetilde{C}, {\cal I}_{p_1,p_2}\E_{|\widetilde{C}})=H^1(C^\p, {\cal I}_{R^j}\E)=0$,
then the restriction map
$H^0(C^\p,\E)\lra H^0(R^j,\E_{|R^j})$ is surjective, so that, by {\rm i)},
the restriction map $H^0(C^\p,\E)\lra H^0(R_\iota^j,\E_{|R_\iota^j})$
to any component $R^j_\iota$ of $R^j$ is surjective, too.
\end{enumerate}
\end{Prop}
\begin{proof}
Part ii) is Proposition 3 in \cite{Sesh}.
\par
Ad i): The restriction of $\F$ to a component $R_\iota$ is of the form
$$
{\Oh}_{{\Pe}_1}(a_1)
\oplus\cdots\oplus{\Oh}_{{\Pe}_1}(a_{\rk \F_{|R_\iota}})\oplus {\rm Torsion}
$$
with
$a_i\ge 0$, $i=1,...,\rk \F_{|R_\iota}$. Therefore, $H^1(R_\iota,\F_{|R_\iota}(-1))=0$.
By successively removing components which are attached at one point only,
the result becomes an easy induction on the length of $R$.
\end{proof}
\begin{Rem}
In the situation of Proposition~\ref{NSproperties} ii), there are
precise conditions for $\pi_*(\E)$ to be torsion free
(\cite{Sesh}, Proposition~5). In particular, any chain of
rational curves contracted by $\pi$ can have length at most
$\rk\E$. This already bounds the family of semistable curves
which might appear in our investigations.
\end{Rem}
\begin{Prop}
\label{NSemb}
Let $C$ be a semistable curve containing the disjoint chains $R^1,...,R^c$ of projective
lines which are attached at two points only, define $\widetilde{C}_j$ as the closure of $C\setminus R^j$, and
let $p_1^j,p_2^j$ be the points where $R^j$ is attached, $j=1,...,c$.
Also set $R:=\bigcup_{j=1}^c R^j$ and define $\widetilde{C}$
as the closure of $C\setminus R$.
Suppose $\E$ is a vector bundle on $C$ which satisfies the following properties
\begin{enumerate}
\item The restriction $\E_{|R^j_\iota}$ of $\E$ to any component of $R^j$ has positive degree,
      $j=1,...,c$.
\item $H^1(\widetilde{C}_j, {\cal I}_{p_1^j,p_2^j}\E_{|\widetilde{C}_j})=0$,
      $j=1,...,c$.
\item The homomorphism $H^0(\widetilde{C}_j, {\cal I}_{p_1^j,p_2^j}\E_{|\widetilde{C}_j})
      \lra \bigl({\cal I}_{p_1^j,p_2^j}\E_{|\widetilde{C}_j}\bigr)/
      \bigl({\cal I}^2_{p_1^j,p_2^j}\E_{|\widetilde{C}_j}\bigr)$ is surjective, $j=1,...,c$.
\item For any point $x\in \widetilde{C}\setminus \{\, p_1^j,p_2^j,\ j=1,...,c\,\}$, the homomorphism
$$
H^0\bigl(C, {\cal I}_R\E\bigr)\lra \E_{|\widetilde{C}}\big/\Bigl({\cal I}_x^2\E_{|\widetilde{C}}\Bigr)
$$
is surjective. 
\item For any two points $x_1\neq x_2\in \widetilde{C}\setminus \{\, p_1^j,p_2^j,\ j=1,...,c\,\}$, 
the homomorphism
$$
H^0\bigl(C, {\cal I}_R\E\bigr)\lra \E_{|\{x_1\}}\oplus \E_{|\{x_2\}}
$$
is surjective. 
\end{enumerate}
Then, the evaluation map ${\rm ev}\colon H^0(C,\E)\otimes{\Oh}_C\lra\E$ yields a closed
embedding
$$
C\hookrightarrow {\rm Gr}\bigl(H^0(C,\E),\rk\E\bigr).
$$
\end{Prop}
\begin{proof} This is proved like Proposition~4 in \cite{Sesh}.
\end{proof}
\begin{Rem}
Note that the conditions (2) - (5) will be satisfied for $\E(l)$, $l\gg 0.$
\end{Rem}
\begin{Prop}[\cite{Sesh}, Lemma 4]
\label{Rest}
Suppose we are given a commutative diagram
$$
\begin{CD}
Z @> \pi >> W
\\
@VpVV @VVqV
\\
T @= T
\end{CD}
$$
in which $p$ and $q$ are projective, $p$ is flat, and $\pi_*{\Oh}_Z={\Oh}_W$. Assume that
$\E_Z$ is a vector bundle on $Z$, such that, for every point $t\in T$, we
have $R^i\pi_{t*}(\E_{Z|p^{-1}(t)})=0$ for all $i>0$, $\pi_t\colon p^{-1}(t)\lra
q^{-1}(t)$ being the induced morphism. Then,
$$
\pi_*(\E_Z)_{|q^{-1}(t)}\q=\q \pi_{t*}(\E_{Z|p^{-1}(t)}),\q \hbox{for all closed points
$t\in T$}.
$$
\end{Prop}
\begin{Rem}
\label{flat}
Let ${\Oh}_W(1)$ be a $q$-ample line bundle. Assume
$$
H^i\bigl(p^{-1}(t), (\E_{Z}
\otimes\pi^*{\Oh}_W(n))_{|p^{-1}(t)}\bigr)=0,\q  \hbox{
for $n\gg 0$, $i>0$, and all $t\in T$}.
$$
Then, $q_*\bigl((\pi_*\E_Z)(n)\bigr)=p_*\bigl(\E_{Z}
\otimes\pi^*{\Oh}_W(n)\bigr)$ is locally free for all $n\gg 0$, whence
$\pi_*\E_Z$ is $T$-flat.
\end{Rem}
\subsection{Modules over Discrete Valuation Rings}
Let $(R, v\colon R\lra\Z)$ be a discrete valuation ring with uniformizing parameter
$t$, i.e., $v(t)=1$. A finitely generated module $M$ over $R$ will be called
\it almost torsion free\rm, if its torsion submodule is annihilated by $t$.
Likewise, a submodule $M^\p$ of a free module $M$ of finite rank is called
\it almost saturated\rm, if the module $M/M^\p$ is almost torsion free.
\begin{Prop}
\label{diagonal}
Let $M=R^{\oplus r}$ be a free module of rank $r$ and
$$
(*)\qquad 0\subsetneq M_1\subsetneq M_2\subsetneq\cdots\subsetneq M_n\subsetneq M
$$
a chain of almost saturated submodules. Then, there is a basis
$e_1,...,e_r$ for $M$, such that $M_i$ is generated by
$\rho^i_1 e_1,...,\rho^i_{\rk M_i} e_{\rk M_i}$ for appropriate elements
$\rho^i_j\in R$ with $v(\rho_j^i)\in\{\,0,1\,\}$, $j=1,...,\rk M_i$, $i=1,...,n$.
\end{Prop}
Note that we do not require $\rk(M_i)>\rk(M_{i-1})$, $i=2,...,n$.
\begin{proof}
We carry out an induction over $r$, the case $r=1$ being clear.
For the first submodule $M_1$, we may find a basis $e_1^*,...,e_r^*$
with the asserted property \cite{Sh}, 10.5. Take $e_1:=e_1^*$.
Then,
$$
(**)\qquad 0\subsetneq \widetilde{M}_1:=M_1/(M_1\cap\langle e_1\rangle)\subsetneq\cdots\subsetneq
\widetilde{M}_n:=M_n/(M_n\cap\langle e_1\rangle)\subsetneq \widetilde{M}:=M/\langle e_1\rangle
$$
is a filtration of the free module $M/\langle e_1\rangle$ of rank $r-1$.
We claim that $\widetilde{M}_i$ is an almost saturated submodule of
$\widetilde{M}$, $i=1,...,n$. If $e_1\in M_i$, this is just the isomorphism theorem.
Otherwise, $te_1\in M_i$ and $M_i^\p:=\widetilde{M}/\widetilde{M}_i$ is a quotient of
$M_i^{\p\p}:=M/M_i=(M/\langle te_1\rangle)/\widetilde{M}_i$. As $M_i^\p$
and $M_i^{\p\p}$ have the same rank, no element of the free part of
$M_i^{\p\p}$ can map to the torsion of $M_i^{\p}$, i.e., the torsion
of $M_i^{\p}$ is a quotient of the torsion of $M_i^{\p\p}$ which implies the claim.
\par
Therefore, we can apply the induction hypothesis to $(**)$. Let $e_2^\p,...,e_r^\p$
be any lift of the appropriate basis $\ol{e}_2,...,\ol{e}_r$
for $\widetilde{M}$. Suppose that we have
already found a basis of the form $e_1,...,e_{\rk M_i}, e^\p_{\rk M_i+1},...,e_r^\p$,
so that the assertion holds for $M_1,...,M_i$ and $e_2,...,e_{\rk M_i}$ also
lift $\ol{e}_2,...,\ol{e}_{\rk M_i}$. Then, $M_{i+1}$ is spanned w.r.t.\ that basis
by vectors of the form
$$
v_j\q=\q\left(\begin{array}{c} *_j\\ 0\\ \vdots\\ 0 \\ q_j \\ 0 \\ \vdots \\ 0
\end{array}\right)
$$
with $q_j,*_j\in R$, $v(q_j)\in \{\,0,1\,\}$, and $q_j$ at the $j$-th place,
$j=1,...,\rk M_{i+1}$.
If $v(q_1)\le v(*_j)$, we may clearly set $*_j=0$.
If $v(q_1)> v(*_j)$, i.e., $v(q_1)=1$ and $v(*_j)=0$, and also $v(q_j)=0$,
we set $e_j:= v_j$. In fact, this case can only occur, if $j>\rk M_i$ or in all the
$M_\iota$ with $\iota\le i$, $t\cdot e_j$ was the respective basis vector.
But then, the result still holds for $M_1,...,M_i$ and the basis
$e_1,..., e_{j-1}, v_j, e_{j+1},...,e_{\rk M_i}, e^\p_{\rk M_i+1},...,e_r^\p$.
Finally, the case $v(q_1)=1$, $v(*_j)=0$, and $v(q_j)=1$ cannot occur.
Indeed, in that case, the class of the vector
$(0,...,t,...,0)^T$, $t$ at the $j$-th place, in $M/M_{i+1}$ is non-zero, i.e.,
the class of $(0,...,1,...,0)^T$ in $M/M_{i+1}$ is not annihilated by $t$, a contradiction.
\end{proof}
\section{Proof of the Main Theorem}
\subsection{Construction of the Hilbert compactification}
The set of pairs $(\nu\colon\widetilde{C}\lra C,\nu^*\E)$, $\nu$ being a partial
normalization of the stable curve $C$ and $\E$ a semistable sheaf of rank $r$ with
Euler characteristic $\chi$ is clearly bounded. Therefore, we may find an $l_1$, such that,
for every $l\ge l_1$, the following assumptions are met.
\begin{Ass}
\label{genasi0}
For every stable curve $C$, every partial normalization $\nu\colon\widetilde{C}\lra C$,
resolving the nodes, say, $N_1,...,N_\nu$,
every semistable sheaf $\E$ on $C$, any two points
$p_1,p_2\in \widetilde{C}$ mapping to a node of $C$, and $Z:=\nu^{-1}\bigl(\{\, N_1,...,N_\nu\,\}\bigr)$, one has
\begin{enumerate}
\item $H^1(\widetilde{C}, {\cal I}_{p_1,p_2}\nu^*\E(l))=0$.
\item The homomorphism $H^0\bigl(\widetilde{C}, {\cal I}_{p_1,p_2}\nu^*\E(l)\bigr)
      \lra \bigl({\cal I}_{p_1,p_2}\nu^*\E(l)\bigr)/\bigl({\cal I}^2_{p_1,p_2}\nu^*\E(l)\bigr)$
      is surjective.
\item For every $x\in \nu^{-1}\bigl(C\setminus\{\,N_1,...,N_\nu\,\}\bigr)$,
      the homomorphism 
	  $$
	  H^0\bigl(\widetilde{C}, {\cal I}_{Z}\nu^*\E(l)\bigr)
      \lra \nu^*\E(l)/\bigl({\cal I}^2_{x}\nu^*\E(l)\bigr)
	  $$
	  is surjective.
\item For any two points $x_1\neq x_2\in \nu^{-1}(C\setminus\{\,N_1,...,N_\nu\,\})$,
      the homomorphism
      $$H^0\bigl(\widetilde{C}, {\cal I}_{Z}\nu^*\E(l)\bigr)
      \lra \nu^*\E(l)_{|\{x_1\}}\oplus \nu^*\E(l)_{|\{x_2\}}$$
      is surjective.
\end{enumerate}
\end{Ass}
In the following, $l$ is assumed to be at least $l_1$.
Let ${\frak H}:={\frak H}(g;\chi_l,r)$ be as in the introduction, and let
${\cal C}_{\frak H}\hookrightarrow {\frak H}\times {\frak G}$ be the universal
curve. Let $\widetilde{\frak H}$ be the open subset of points $h$ for which
$C_h:={\cal C}_{\frak H}\times_{\frak H} \{h\}$ is semistable.
Let $q_{\widetilde{\frak H}}\colon V^{\chi_l}\otimes{\Oh}_{\cal C_{\widetilde{\frak H}}}
\lra \E_{\widetilde{\frak H}}$ be the pullback of the universal quotient.
Then, there is a flat family $q\colon \cal C^*_{\widetilde{\frak H}}\lra {\widetilde{\frak H}}$
of stable
curves together with an ${\widetilde{\frak H}}$-morphism $\pi\colon
\cal C_{\widetilde{\frak H}}\lra \cal C^*_{\widetilde{\frak H}}$,
such that $\pi_*{\Oh}_{\cal C_{\widetilde{\frak H}}}={\Oh}_{\cal C^*_{\widetilde{\frak H}}}$
and $\pi$ is fibrewise the contraction onto the stable model.
By \ref{NSproperties} ii) (2),
we are in the position to apply Proposition~\ref{Rest}. Moreover, we see
that the assumptions of Remark~\ref{flat} are satisfied.
We get the homomorphism
$$
\pi_*(q_{\widetilde{\frak H}})\colon V^{\chi_l}\otimes{\Oh}_{\cal C^*_{\widetilde{\frak H}}}
\lra \pi_*\E_{{\widetilde{\frak H}}}
$$
of ${\widetilde{\frak H}}$-flat sheaves. Let ${\frak H}^0$ be the open set of points
$h$ for which
\begin{itemize}
\item $\bigl(\pi_*\E_{\widetilde{\frak H}}\bigr)_{|q^{-1}(h)}$ is a semistable
      sheaf.
\item $\pi_*(q_{\widetilde{\frak H}})_{|q^{-1}(h)}$ is surjective.
\item $H^0(\pi_*(q_{\widetilde{\frak H}})_{|q^{-1}(h)})$ is an isomorphism.
\end{itemize}
Set ${\cal C}^*_{\frak H^0}:=
{\cal C}^{*}_{\widetilde{\frak H}}\times_{\widetilde{\frak H}}{\frak H^0}$.
Then, the quotient family
$$
\pi_{\frak H^0}:=
\pi_*(q_{\widetilde{\frak H}})_{|\cal C^*_{\frak H^0}}
\colon V^{\chi_l}\otimes{\Oh}_{\cal C^*_{{\frak H}^0}}
\lra \E_{{\frak H^0}}:=(\pi_*\E_{{\widetilde{\frak H}}})_{|\cal C^*_{\frak H^0}}
$$
defines a morphism ${\frak H}^0\lra \ol{\frak Q}(g;\chi_l,r)$.
Let ${\frak X}\subset {\frak H}\times \ol{\frak Q}(g;\chi_l,r)$ be the closure
of the graph of the above morphism. Let ${\frak L}_{a,m}$ be as in the introduction
with $a$ so large, that the semistable points in ${\frak X}$ lie in the preimage
of the points in $\ol{\frak Q}(g;\chi_l,r)$ which are semistable w.r.t.\ the linearization
in ${\frak N}$.
\begin{Thm}
\label{Semstabred}
Let $x\in {\frak X}$ be a point which is semistable w.r.t.\ the linearization
in ${\frak L}_{a,m}$, $a\gg 0$. Then, $x$ is contained in the graph $\Gamma$
of the morphism ${\frak H}^0\lra \ol{\frak Q}(g;\chi_l,r)$.
\end{Thm}
\begin{proof}
By construction of ${\frak H}$ and ${\frak X}$, the set of
(semi)stable points $y\in {\frak X}$ corresponding to smooth
curves is dense. Let $y$ represent the smooth curve $C$. By the
result of \cite{Sch}, $y$ is (semi)stable (w.r.t.\ the
linearization in ${\frak L}_{a,m}$ for all $a\gg 0$), if and only
if $V^{\chi_l}\lra H^0(\E(l))$ is an isomorphism and $\E$ is a
(semi)stable bundle of rank $r$ with Euler characteristic $\chi$,
$\E:=\E_{{\frak H}^0|C}$. Hence, we may find a smooth curve $K$,
a point $k\in K$, and a morphism $\kappa\colon K\lra {\frak X}$
with $\kappa(k)=x$, such that $\kappa(K\setminus\{k\})$ is
contained in the locus of pairs $(C,q_C\colon
V^{\chi_l}\otimes{\Oh}_C\lra\E(l))\in {\Gamma}$ for which $C$ is a
smooth curve and $\E$ is a stable vector bundle. We have an
induced morphism $\ol{\kappa}\colon K\lra \ol{\frak
Q}(g;\chi_l,r)$ which lands, by assumption, in the semistable
locus. Without loss of generality, we may assume that this
morphism is induced by a family $(\cal C_K^*, q_K\colon
V^{\chi_l}\otimes{\Oh}_{\cal C^*_K}\lra\E_{\cal C^*_K})$. The
surface $\sigma\colon S:=\cal C_K^*\lra K$ is smooth outside the
nodes of $\sigma^{-1}(k)$ and has singularities of type $A_n$ in
these nodes. We may resolve these singularities in the usual way
in order to get a flat family $\widetilde{\sigma}\colon
\widetilde{S}=\widetilde{\cal C}_K\lra K$ of semistable curves
with $\widetilde{S}$ smooth. Let $\widetilde{q}_K$ be the
pullback of $q_K$ to $\widetilde{S}$. Then, $\widetilde{q}_K$
defines a rational map $\widetilde{S}\lra {\frak G}$ which is
defined outside some nodes of $\widetilde{\sigma}^{-1}(k)$. By
blowing up these nodes and points which are infinitely near to
them, we get a new flat family $\sigma^\p\colon S^\p=\cal
C^\p_K\lra K$ and $q^\p_K\colon V^{\chi_l}\otimes{\Oh}_{S^\p}
\lra\E^\p_K$ where $\E_K^\p$ is locally free. Set
$C_k^\p:={\sigma^\p}^{-1}(k)$. Let $\widetilde{C}$ be the closure of
$C^\p_k$ with all chains of rational curves attached at only two
points removed. Then, $\nu\colon \widetilde{C}\lra C$ is a
partial normalization of $C$. By construction, the morphisms
$\widetilde{C}\lra {\frak G}$ induced by $q^\p_{K|\widetilde{C}}$
and $\nu^*(q_{K|\sigma^{-1}(k)})$ agree, whence these quotients
are equivalent. Now, our Assumptions~\ref{genasi0} and
Proposition~\ref{NSemb} imply that the image of $C_k^\p$ under
the map $S^\p\lra {\frak G}$ is a semistable curve $C_k^{\p\p}$,
and $\pi\colon C_k^\p\lra C_k^{\p\p}$ just contracts all rational
chains on which $\E_K^\p$ is trivial. Next, look at the
commutative diagram
$$
\begin{CD}
\widetilde{C} @=  \widetilde{C} @>\nu>> C
\\
@| @V\subset VV @|
\\
\widetilde{C} @>\subset >> C_k^{\p} @> \nu^\p >> C
\\
@| @V\pi VV @|
\\
\widetilde{C} @>\subset >> C_k^{\p\p} @>\nu^{\p\p}>> C.
\end{CD}
$$
The quotients
$\pi_*(q^\p_{K|{\sigma^\p}^{-1}(k)}\bigr)_{|\widetilde{C}}$ and $\nu^*(q_{K|\sigma^{-1}(k)})$ are also equivalent.
Thus, there is an isomorphism
$\alpha\colon \bigl(\pi_*\E^\p_{K|{\sigma^\p}^{-1}(k)})_{|\widetilde{C}}
\lra \nu^*\bigl(\E_{K|\sigma^{-1}(k)}\bigr)$, making the following diagram commute
$$
\begin{CD}
V^{\chi_l}\otimes{\Oh}_{C^{\p\p}_k} @>\q \q\qquad
\pi_*(q^\p_{K|{\sigma^\p}^{-1}(k)})\qquad\q\q >>
\pi_*\E^\p_{K|{\sigma^\p}^{-1}(k)}
\\
@V\hbox{restriction to $\widetilde{C}$}VV @VV\hbox{restriction to $\widetilde{C}$}V
\\
V^{\chi_l}\otimes{\Oh}_{\widetilde{C}}
@>\q\qquad \pi_*(q^\p_{K|{\sigma^\p}^{-1}(k)})_{|\widetilde{C}}\qquad\q
>>  \bigl(\pi_*\E^\p_{K|{\sigma^\p}^{-1}(k)}\bigr)_{|\widetilde{C}}
\\
@| @VV\alpha V
\\
V^{\chi_l}\otimes{\Oh}_{\widetilde{C}}
@>\q\q\q\qquad \nu^*(q_{K|\sigma^{-1}(k)})\qquad\q\q\q
>>  \nu^*\E_{K|{\sigma}^{-1}(k)}.
\end{CD}
$$
This latter diagram finally provides us, via projection onto $C$,
with the following commutative diagram
$$
\begin{CD}
V^{\chi_l}\otimes{\Oh}_C @>\q \q\qquad\nu^{\p\p}_*\pi_*(q^\p_{K|{\sigma^\p}^{-1}(k)})\qquad\q\q >>
\nu^{\p\p}_*\pi_*\E^\p_{K|{\sigma^\p}^{-1}(k)}
\\
@VVV @VVV
\\
V^{\chi_l}\otimes\nu_*{\Oh}_{\widetilde{C}}
@> \q\q\qquad\nu_*\nu^*(q_{K|\sigma^{-1}(k)})\q\q\qquad
>>  \E_{K|\sigma^{-1}(k)}\otimes\nu_*{\Oh}_{\widetilde{C}}.
\end{CD}
$$
Furthermore, we have the commutative diagram
$$
\begin{CD}
V^{\chi_l}\otimes{\Oh}_C @>\q q_{K|\sigma^{-1}(k)}\q>> \E_{K|\sigma^{-1}(k)}
\\
@VVV @VVV
\\
V^{\chi_l}\otimes\nu_*{\Oh}_{\widetilde{C}}
@> \nu_*\nu^*(q_{K|\sigma^{-1}(k)})>>  \E_{K|\sigma^{-1}(k)}\otimes\nu_*{\Oh}_{\widetilde{C}}
\end{CD}
$$
in which the vertical arrows are injective.
Therefore, the image of $\nu^{\p\p}_*\pi_*\E^\p_{K|{\sigma^\p}^{-1}(k)}$
in the sheaf $\E_{K|\sigma^{-1}(k)}\otimes\nu_*{\Oh}_{\widetilde{C}}$ is
$\E_{K|\sigma^{-1}(k)}$. The kernel of the surjection
$$
\nu^{\p\p}_*\pi_*\E^\p_{K|{\sigma^\p}^{-1}(k)}\lra \E_{K|\sigma^{-1}(k)}
$$
must be zero, because both sheaves have the same Hilbert polynomial w.r.t.\ ${\Oh}_C(1)$.
\par
To conclude, the quotient $q_K^\p$ defines a $K$-morphism
$\vartheta\colon {\cal C}^\p_K\lra K\times {\frak G}$ and we have
seen (a) that the image is a flat family of curves
$\sigma^{\p\p}\colon {\cal C}^{\p\p}_K\lra K$ with $C^{\p\p}$ as
the fibre over $k$ and (b) that the family
$q^{\p\p}_K:=\vartheta_*(q_K^\p)\colon V^{\chi_l}\otimes
{\Oh}_{{\cal C}_K^{\p\p}} \lra \E_{{\cal
C}_K^{\p\p}}:=\vartheta_*\E_{{\cal C}_K^\p}$ is a flat quotient
(by Proposition~\ref{Rest}), such that
$\sigma^{\p\p}_*(q_K^{\p\p})$ is a quotient onto a family of
semistable torsion free coherent sheaves. Therefore, the family
${\cal C}_K^{\p\p}$ defines a morphism $\kappa^\p\colon K\lra
\Gamma\subset {\frak X}$. Since we have not altered the original
family outside the point $k$, $\kappa^\p$ agrees with the
original morphism $\kappa$ outside $k$, and, thus, everywhere,
because ${\frak X}$ is separated. This shows $x=\kappa^\p(k)\in
\Gamma$.
\end{proof}
\begin{Rem}
One could also use the arguments presented in the paper
\cite{Sesh}. We have chosen the alternative way, also suggested
in \cite{Sesh}, because it reflects more of the moduli problem we
are dealing with.
\end{Rem}
Thus, the GIT-quotient ${\frak H}^0\catqot_{\frak L_{a,m}}\SL(V^{\chi_l})$ exists as a
projective scheme over $\ol{\frak M}_g$. It also comes with an $\ol{\frak M}_g$-morphism
to $\ol{\frak U}(g;\chi,r)$. The hard task will be to give it a modular interpretation.
\subsection{Analysis of semistability}
Fix the data $g$, $\chi$, and $r$, and set $\chi_l:=\chi+l\cdot
r\cdot (2g-2)$. Let $C$ be a semistable curve of genus $g$.
Denote by $\pi\colon C\lra C^\p$ the contraction onto the stable
model of $C$. Let $\E$ be a vector bundle of uniform rank $r$ on
$C$ with Euler characteristic $\chi(\E) =\chi$, such that $\E$
has positive degree on each rational component and $\pi_*(\E)$ is
torsion free. From Proposition~\ref{NSemb}, we infer that, for
sufficiently large $l$, $H^0(\E(l))\otimes{\Oh}_C\lra \E(l)$ will
give rise to a closed embedding $C\hookrightarrow {\frak
G}(H^0(\E(l)),r)$ and $H^1(\E(l))=\{0\}$. Identifying $H^0(\E(l))$
with some fixed vector space $V^{\chi_l}$ of dimension $\chi_l$,
we may ask whether $[C]\in {\frak H}^0(g;\chi_l,r)$ is
$\SL(V^{\chi_l})$-(semi)stable w.r.t.\ the linearization in
${\frak L}_{a,m}$ for large $a$. We already know from
Pandharipande's construction and \ref{suitable} that (a) if
$\pi_*(\E)$ is stable w.r.t.\ the canonical polarization, then
$[C]$ will be stable and (b) if $[C]$ is semistable, then
$\pi_*(\E)$ is semistable w.r.t.\ the canonical polarization. If
$\pi_*(\E)$ is properly semistable, then there will be additional
conditions for $[C]$ to be (semi)stable. We will have to analyze
those conditions. Abstractly, by \ref{Absi}, they can be
described as follows: Suppose we are given  $q_C\colon
V^{\chi_l}\otimes {\Oh}_C\lra\E(l)$, such that $H^0(q_C)$ is an
isomorphism and $\pi_*(\E)$ is semistable w.r.t.\ the canonical
polarization. Then, $[C]$ will be $\SL(V^{\chi_l})$-(semi)stable,
if and only if for every one parameter subgroup
$\la\colon\C^*\lra \SL(V^{\chi_l})$ with $\mu_{\frak N}(\la,
[\pi_*(q_C)])=0$, one has
\begin{equation}
\label{condition}
\mu_{\frak L_m}\bigl(\la, \wedge^{P(m)}\Psi_{\frak H(g;\chi_l,r)}\bigr)
\q (\ge)\q 0,\qquad  \Psi_{\frak H(g;\chi_l,r)}\hbox{ as in the introduction}.
\end{equation}
Let us first recall when $\mu_{\frak N}(\la,[\pi_*(q_C)])=0$ happens.
For this, suppose $\la$ is given with respect to the basis $v_1,...,v_{\chi_l}$ by the
weight vector
\begin{equation}
\label{weight}
\ul{\gamma}\q=\q \sum_{i=1}^{{\chi_l}-1}\alpha_i\bigl(\underbrace{i-{\chi_l},...,i-{\chi_l}}_{i\times},\underbrace
{i,...,i}_{{\chi_l}\times}\bigr),\q \alpha_i\in\Q_{\ge 0},\ i=1,...,{\chi_l}-1,
\end{equation}
and let $i_1<i_2<\cdots <i_k$ be the indices with $\alpha_{i}>0$.
Then, we get a filtration $V^\bullet:\{0\}=:V_0\subset
V_1\subset\cdots\subset V_k\subset V_{k+1}:=V^{\chi_l}$ with
$V_j:=\langle\, v_1,...,v_{i_j}\,\rangle$, $j=1,...,k$. Recall
that $\mu_{\frak N}(\la, [\pi_*(q_C)])$ depends only on the pair
$(V^\bullet,\ul{\alpha})$,
$\ul{\alpha}:=(\alpha_{i_1},...,\alpha_{i_k})$ (\cite{GIT},
Section~2.2). Set
$$
\widetilde{\F}_j\q:=\q\pi_*\Bigl(q_C\bigl(V_j\otimes {\Oh}_C\bigr)(-l)\Bigr),\qquad j=1,...,k.
$$
Now, from the construction of $\ol{\frak U}(g;\chi_l,r)$ and Lemma~\ref{equal},
one knows that the equality
$\mu_{\frak N}(\la,[\pi_*(q_C)])=0$ will occur, if and only if
\begin{itemize}
\item $\sum_{j=1}^k\alpha_{i_j}\bigl(P(\pi_*(\E))\trk(\widetilde{\F}_j)-P(\widetilde{\F}_j)\trk(\E)\bigr)=0$,
      i.e., each $\widetilde{\F}_j$ destabilizes $\E$;
\item $H^0(\pi_*(q_C))(V_j)=H^0(\widetilde{\F}_j(l))$, $j=1,...,k$.
\end{itemize}
Recall that the family of pairs $(C,\F)$ with $C$ a stable curve and $\F$ a destabilizing
subsheaf of a semistable torsion free sheaf $\E$ of uniform rank $r$ on $C$ with
$\chi(\E)=r$ is bounded. In view of \ref{NSproperties}, we can find an $l_1$, such that
for all $l\ge l_1$, the following assumptions are verified.
\begin{Ass}
\label{genasi}
Let $[C]\in {\frak H}^0(g;\chi_l,r)$ be semistable w.r.t.\ the linearization
in ${\frak L}_{a,m}(g;\chi_l,r)$ for $a\gg 0$. Denote by $\pi\colon C\lra C^\p$
the contraction onto the stable model, and by $q_C\colon V^{\chi_l}\otimes{\Oh}_C\lra
\E$ the pullback of the universal quotient. Then, $\pi_*(\E(-l))$ is, by definition,
semi\-stable and $H^0(\pi_*(q_C))$ is an isomorphism.
For every destabilizing subsheaf $\F^\p\subset\E(-l)$, define
$$
\F\q:=\q q_C\biggl(
\Bigl(H^0\bigl(\pi_*(q_C)\bigr)^{-1}\bigl(H^0(\F^\p(l))\bigr)\Bigr)\otimes{\Oh}_C\biggr).
$$
Note that $\pi_*(\F(-l))=\F^\p$.
Then, we assume:
\begin{enumerate}
\item For every irreducible component $C_\iota$, the restriction map
      $H^0(C,\F)\lra H^0(C_{\iota},\F_{|C_\iota})$ is surjective.
\item For every irreducible component $C_\iota$ which is not a rational curve attached at only two points
      and on which $\F$ has positive rank, and
      any two points $p_1\neq p_2\in C_\iota$, the evaluation map
      $H^0(C_{\iota},\F_{|C_\iota})\lra \F_{|C_\iota}\otimes_{{\Oh}_C}{\Oh}_{\{\,p_1,p_2\,\}}$
      is surjective.
\item For every irreducible component $C_\iota$ which is not a rational curve attached at only two points
      and on which $\F$ has positive rank and
      any point $p\in C_\iota$, the evaluation map
      $H^0(C_{\iota},\F_{|C_\iota})\lra \F_{|C_\iota}\otimes_{{\Oh}_C}\bigl({\frak m}_{p,C_\iota}
      /{\frak m}_{p,C_\iota}^2\bigr)$ is surjective.
\end{enumerate}
\end{Ass}
Next, we fix a maximal filtration $\F_\bullet\colon
0=:\F_0^\p\subsetneq \F_1^\p\subsetneq\F_2^\p\subsetneq
\cdots\subsetneq\F_k^\p\subsetneq\F_{k+1}^\p:=\pi_*(\E)$ of
$\pi_*(\E)$ by destabilizing subsheaves and a vector
$\ul{\alpha}:=(\alpha_1,...,\alpha_k)$ of non-negative rational
numbers. Let us also fix some $l$ and define
$V_j:=H^0(\F^\p_j(l))$ under the identification of $V^{\chi_l}$
with $H^0(\pi_*(\E)(l))$ and $\F_j:=q_C(V_j\otimes{\Oh}_C)$,
$j=1,...,k$.
\begin{Rem}
\label{satt}
The $\F_j$ are saturated subsheaves of $\E$, $j=1,...,k$. For, if $\widetilde{\F}_j$
is the saturation of $\F_j$, then $\widetilde{\F}_j(l^\p)$ will be globally generated for
some $l^\p$ large enough. This is clear for points which lie on a component on
which $\omega_C$ is ample. Now, let $R$ be the disjoint union of all maximal
chains of rational curves attached at only two points, set
$C^*:=C\setminus R$, let $\ol{C}^*$
be the closure of $C^*$, and define $\ul{x}:=\ol{C}^*\setminus C^*$.
Then, $H^1\bigr(\ol{C}^*, {\cal I}_{\ul{x}}
\widetilde{\F}_j(l^\p)\bigr)$ is zero for $l^\p\gg 0$, i.e.,
$H^0\bigl(C,\widetilde{\F}_j(l^\p)\bigr)\lra H^0\bigl(R, \widetilde{\F}_{j|R}(l^\p)\bigr)$
is surjective.
As $\widetilde{\F}_{j|R}(l^\p)$ is globally generated, the claim is settled.
Now, if $\widetilde{\F}_j$ strictly contains $\F_j$, then
$h^0\bigl(C,\widetilde{\F}_j(l^\p)\bigr)>h^0\bigl(C,\F_j(l^\p)\bigr)$ for all $l^\p\gg 0$.
This is, however, not possible,
because the inclusion $\F^\p_j(l+l^\p)=\pi_*\F_j(l^\p)\subset \pi_*\widetilde{\F}_j(l^\p)$
must be an equality as both sheaves have the same multi-rank and $\F^\p_j$ is saturated.
\end{Rem}
Note that, given a basis
$v_1,...,v_{\chi_l}$ of $V^{\chi_l}$ such that
$\langle\, v_1,...,v_{\dim V_j}\,\rangle=V_j$, $j=1,...,k$,
the vector $\ul{\alpha}$ defines
a one parameter subgroup $\la$ of $\SL(V^{\chi_l})$, by Formula (\ref{weight}).
By \cite{GIT}, 2.2,  again, the values $\mu_{\frak N}(\la,[\pi_*(q_C)])$ and
$\mu_{\frak L_m}(\la, [\wedge^{P(m)}\Psi_{\frak H(g;\chi_l,r)}])$ do not depend on the
choice of such a basis. It will be our task to compute the value of
$\mu_{\frak L_m}(\la, [\wedge^{P(m)}\Psi_{\frak H(g;\chi_l,r)}])$.
In other words, we will have to find a basis for
$H^0(C,\det(\E(l)^{\otimes m}))$ for which the weight of the associated element in
$\Hom(\wedge^{P(m)} S^m V,\C)$ becomes minimal.
Then, $\mu_{\frak L_m}(\la,\wedge^{P(m)}\psi)$ will be minus that weight.
\subsubsection*{Case {\rm A)}: The irreducible components of $C$ are smooth}
Write $C=\bigcup_{\iota=1}^c C_\iota$ as the union of its irreducible components
and set $\L_m:=\det(\E(l))^{\otimes m}$.
For large $m$, the restriction map $H^0(C,\L_m)\lra H^0(C_\iota,\L^\iota_m)$, $\L^\iota_m:=
\L_{m|C_\iota}$, will be surjective for $\iota=1,...,c$.
We may, therefore, compute first the weights of a basis for $H^0(C_\iota,\L_m^\iota)$,
$\iota=1,...,c$.
Set $\E_\iota:=\E(l)_{|C_\iota}$ and $\F_j^\iota:=\Im(\F_j\lra \E_\iota)$, $j=1,...,k$,
$\iota=1,...,c$.
We also define $\widetilde{V}_j:=\langle\, v_{\dim V_{j-1}+1},...,v_{\dim V_j}
\,\rangle$, $j=2,...,k+1$,
and $\widetilde{V}_1:=V_1$.
Then,
$$
\bigwedge^r V=\bigoplus_{(\rho_1,...,\rho_{k+1}): \sum \rho_j=r} W_{\rho_1,...,\rho_{k+1}},
\q W_{\rho_1,...,\rho_{k+1}}:=\bigwedge^{\rho_1} \widetilde{V}_1
\otimes\bigwedge^{\rho_2}\widetilde{V}_2
\otimes\cdots\otimes\bigwedge^{\rho_{k+1}}\widetilde{V}_{k+1}.
$$
The spaces $W_{\ul{\rho}}$ are weight spaces for $\la$ for the
weight
$$
w^{\ul{\alpha}}_{\iota,\ul{\rho}}(l)\q:=\q {\rho_{1}}\cdot\gamma_1(l)+\cdots+{\rho_{k+1}}\cdot\gamma_{k+1}(l)
$$
where $\gamma_j(l)$ is the weight of a section in $\widetilde{V}_j$, $j=1,...,k+1$.
Formula~(\ref{weight}) shows that the $\gamma_j(l)$ are, in fact, polynomials
in $l$.
\begin{Rem}
Let $N\in C$ be a node of $C$, i.e., a point where two components $C_1$ and $C_2$ of $C$ meet.
Then, the stalk of a torsion free sheaf of ${\cal Q}$ at $N$ is of the form
${\Oh}_{C,N}^{\oplus a_1}\oplus {\Oh}_{C_1,N}^{\oplus a_2}\oplus {\Oh}_{C_2,N}^{\oplus a_3}$
(\cite{Sesh2}, Huiti\^eme Partie, Proposition 3).
Thus, if $\E$ is a torsion free sheaf and $\F$ is a saturated subsheaf, then
the image $\F_1$of $\F$ in $\E_{|C_1}$ is an almost saturated subsheaf, because
${\rm Tors}(\E_{|C_1}/\F_1)\cong \C^{\oplus a_3}$ as ${\Oh}_{C_1,N}$-module.
Therefore, the Structure Result~\ref{diagonal} allows to determine the vanishing orders
of sections coming from a weight space $W_{\ul{\rho}}$. This observation will be
crucial for the following subtle analysis of weights and vanishing orders.
\end{Rem}
Let us look at some specific $\iota\in\{\, 1,...,c\,\}$. Then, the space of minimal weight
which produces sections which do not vanish on $C_\iota$ is
$W^\iota_{\rm min}:=W_{\rk \F_1^\iota,\rk\F_2^\iota-\rk\F_1^\iota,...,\rk\E_\iota-\rk \F_k^
\iota}$.
The associated weight is
\begin{eqnarray*}
w_{\iota,\rm min}^{\ul{\alpha}}(l)&:=&\sum_{j=1}^k\alpha_{j}\bigl(
\dim V_j\cdot \rk \E_\iota-\chi_l\cdot \rk\F_k^\iota \bigr)
\\
&=& \sum_{j=1}^k\alpha_{j}\bigl(P(\F^\p_j)(l)\cdot\rk \E_\iota -
P(\pi_*(\E))(l)\cdot \rk\F_j^\iota\bigr).
\end{eqnarray*}
Let $N^\iota_1,...,N^\iota_{\nu_\iota}$ be the nodes of $C$ located on $C_\iota$, i.e.,
the points where $C_\iota$ meets other components of $C$.
Note that each $\F^\iota_j$ is a subbundle of $\E_\iota$ outside
the above points. This is because $\F_j$ is a saturated subsheaf of $\E$, by Remark~\ref{satt},
$j=1,...,k$. Let us look at a specific node $N\in\{\,N^\iota_1,...,N^\iota_{\nu_\iota}\,\}$.
Let $\e$ be the fibre of $\E$ at $N$, and $\f_j$ the image of $\F^\iota_j$
in $\e$, $j=1,...,k$. Set $a_j:=\dim \f_{j}$, $r_j:=\rk\F^\iota_j$, and
$b_j:=\min\{\, a_j-a_{j-1}, r_j-r_{j-1}\,\}$, $j=1,...,k+1$.
A general section of $W^\iota_{\rm min}$ will vanish of order $o_N:=r-b_N$ at $N$.
This is an immediate consequence of Proposition~\ref{diagonal}.
\begin{Lem}
\label{generate}
{\rm i)}
The sections in $W^\iota_{\rm min}$ generate
$\det\bigl(\E_{|C_\iota}(l)(-\sum_{s=1}^{\nu_{\iota}} o_{N_s} N_s)\bigr)$.
\par
{\rm ii)} The image of $W^\iota_{\rm min}$ in
$H^0\bigl(\det(\E_{|C_\iota}(l)(-\sum_{s=1}^{\nu_{\iota}} o_{N_s} N_s))\bigr)$ is a very ample
linear system, unless $C_\iota$ is rational and attached at only two points
and $\det\bigl(\E_{|C_\iota}(l)(-\sum_{s=1}^{\nu_{\iota}} o_{N_s} N_s)\bigr)$ is trivial.
\end{Lem}
\begin{proof}
The assertion i) about global generation results immediately from Assumption~\ref{genasi} (1).
Likewise, \ref{genasi} (2) and (3) settle the very ampleness in ii) when $C_\iota$ is not a rational
component which is attached at only two points.
In the remaining case, $\E_\iota$ is of the form ${\Oh}_{{\Pe}_1}^{\oplus s}\oplus{\Oh}_{{\Pe}_1}(1)
^{\oplus t}$. The $\F^\iota_j$, as globally generated subsheaves of $\E_\iota$, are also
of the form ${\Oh}_{{\Pe}_1}^{\oplus s_j}\oplus{\Oh}_{{\Pe}_1}(1)^{\oplus t_j}$, $j=1,...,k$.
Let $r_1,...,r_{b+1}$ be the ranks occurring among the $\F^\iota_j$, and for
$\beta=1,...,b+1$, let $\G_\beta$ be the first among the $\F^\iota_j$ to attain that rank.
There are now two possibilities: either $H^0(\G_{\beta+1})/H^0(\G_{\beta})$
contains a subspace $H^0({\Oh}_{{\Pe}_1}(1))$ for some $\beta\in\{\,0,...,b\,\}$ or not.
In the first case, by \ref{genasi} (1), it is easy to explicitly construct sections
separating points and tangent vectors. The second case, however, can only occur,
if all the $\G_\beta$ are trivial, again by \ref{genasi}. But then,
$\det\bigl(\E_{|C_\iota}(l)(-\sum_{s=1}^{\nu_{\iota}} o_{N_s} N_s)\bigr)$ is obviously trivial.
\end{proof}
\begin{Cor}
\label{asi1} The space
$S^m W_{\rm min}^{\iota}$ generates
$H^0\Bigl(\bigl(
\det(\E(l))_{|C_\iota}(-\sum_{s=1}^{\nu_{\iota}} o_{N_s} N_s)\bigr)^{\otimes m}\Bigr)$
for all $m\gg 0$.
\end{Cor}
Let $d^\p_\iota$ be the degree of $\det(\E_{|C_\iota}(-\sum_{s=1}^{\nu_{\iota}} o_{N_s} N_s))$
and $e_\iota:=\deg(\omega_{C|C_\iota})$.
Because of Corollary~\ref{asi1},
the elements in $S^m W_{\rm min}^\iota$ will contribute
the weight
$$
K^{\ul{\alpha}}_\iota(l,m)\q:=\q
m\cdot\bigl(m\cdot(d^\p_\iota+l\cdot r\cdot e_\iota)+ 1-g(C_j)\bigr)\cdot w^{\ul{\alpha}}_
{\iota,\rm min}(l)
$$
to a basis of $H^0(C_\iota, \L_m^\iota)$. Thus, we only have to worry about sections
vanishing of lower order than $m\cdot o_N$ at $N$.
Note that $W_{\rho_1,...,\rho_{k+1}}$ will produce sections which do not vanish on
$C_{\iota}$ if and only if the condition
\begin{equation}
\label{notvanish}
\sum_{i=1}^{j} \rho_i\q\le\q r_{j},\q j=1,...,k,
\end{equation}
is satisfied. A tuple $\ul{\rho}=(\rho_1,...,\rho_{k+1})$ satisfying (\ref{notvanish})
will be called \it admissible\rm.
Note that there are only finitely many admissible tuples.
\par
Next, we let $\kappa_1<\cdots<\kappa_s$ be the elements in $\{\,1,...,k+1\,\}$
with $a_{\kappa}-a_{\kappa-1}>r_\kappa-r_{\kappa-1}$, set $K:=\{\,\kappa_1,...,
\kappa_s\,\}$ and $K^*:=\{\,1,...,k+1\,\}\setminus K$.
\begin{Lem}
\label{weightincrease1}
Fix a vanishing order $o<o_N$, and let $w^{\ul{\alpha}}_{\iota,\ul{\rho}}(l)$
be the minimal weight of a
section with vanishing order $o$. Then, $\ul{\rho}$ may be chosen to satisfy
$$
\begin{array}{cccccl}
r_\kappa-r_{\kappa-1} & \le &\rho_\kappa & \le & a_{\kappa}-a_{\kappa-1}
&\hbox{for $\kappa\in K$}
\\
a_\kappa-a_{\kappa-1} & \le &\rho_\kappa & \le & r_{\kappa}-r_{\kappa-1}
&\hbox{for $\kappa\in K^*$}.
\end{array}
$$
\end{Lem}
\begin{proof}
We begin with the right hand side inequalities. Suppose
$\rho_\kappa$ violates the right hand inequality.
In particular, $\rho_\kappa>r_\kappa-r_{\kappa-1}$, whence
$\sum_{j=1}^{\kappa-1}\rho_j< r_{\kappa-1}$. Define $\ul{\rho}^\p=(\rho_1^\p,...,
\rho_{k+1}^\p)$ by $\rho^\p_{\kappa-1}:=\rho_{\kappa-1}+1$, $\rho_{\kappa}^\p:=
\rho_\kappa-1$, and $\rho_j^\p:=\rho_j$ for $j\neq \kappa-1,\kappa$.
Then, $\ul{\rho}^\p$ is obviously admissible,
$w^{\ul{\alpha}}_{\iota,\ul{\rho}^\p}(l)\le w^{\ul{\alpha}}_{\iota,
\ul{\rho}}(l)$, and $W_{\ul{\rho^\p}}$ will still produce
sections of vanishing order $o$. The latter property results from the fact that $\rho_{\kappa}$ was,
by assumption,
strictly bigger than the maximal number of sections in $\widetilde{V}_\kappa$ with linearly
independent images in $\e$.
\par
The other inequality asserts $\rho_j\ge b_j$ for $j=1,...,k+1$. For $k+1$, we
have $\sum_{j=1}^k\rho_j\le r_k$, so $\rho_{k+1}=r-\sum_{j=1}^k\rho\ge r_{k+1}-r_k$.
Suppose $\rho_{j_0}<b_{j_0}$.
Then, there is an index $j^\p>j_0$ with $\rho_{j^\p}\neq 0$.
Otherwise, $\sum_{j=1}^{j_0}\rho_j=r$ and then $\rho_{j_0}\ge r_{j_0}-r_{j_0-1}$ as before,
a contradiction. Let $j^\p$ be minimal with the above properties. Define
$\ul{\rho}^\p=(\rho_1^\p,...,\rho^\p_{k+1})$ with $\rho^\p_{j_0}:=\rho_{j_0}+1$,
$\rho^\p_{j^\p}:=\rho_{j^\p}-1$, and $\rho^\p_j=\rho_j$ for $j\neq j_0,j^\p$.
As $\sum_{j=1}^{j_0-1}\rho_j\le r_{j_0-1}$ and $\rho_{j_0}<r_{j_0}-r_{j_0-1}$, $\ul{\rho}^\p$
is again admissible. Moreover,
$w^{\ul{\alpha}}_{\iota,\ul{\rho}^\p}(l)\le w^{\ul{\alpha}}_{\iota,\ul{\rho}}(l)$ and $W_{\ul{\rho^\p}}$
will still lead to
sections of vanishing order $o$. This time, the last assertion is the consequence of the
assumption that $\rho_{\kappa}$ was
strictly smaller than the maximal number of sections in $\widetilde{V}_\kappa$ with linearly
independent images in $\e$.
\end{proof}
Let $I\subset K^*\times K$ be the set of all
$(i,j)$ with $i<j$. Fix an order "$\preceq$" on $I$, such that $(i^\p,j^\p)\preceq
(i,j)$ implies $\gamma_{j^\p}(l)-\gamma_{i^\p}(l)\le \gamma_j(l)-\gamma_i(l)$.
The idea for the following investigations is the following: Suppose we are given
an admissible tuple $\ul{\rho}$, satisfying the inequalities of Lemma~\ref{weightincrease1},
and $(i,j)\in I$. Then, we define a new tuple $\ul{\rho}^\p$ with $\rho^\p_i:=\rho_i-1$
and $\rho^\p_j=\rho_j+1$ and let all other entries of $\ul{\rho}^\p$ agree with those
of $\ul{\rho}$. As $i<j$, $\ul{\rho}^\p$ is still admissible. However, we will perform
this operation only if $\ul{\rho}^\p$ still satisfies the inequalities of
Lemma~\ref{weightincrease1}. In that case, the generic vanishing order of sections
in $W_{\ul{\rho}^\p}$ will be one less than the generic vanishing order in $W_{\ul{\rho}}$.
Thus, if we are given a specific vanishing order $o$, we carry out $s:=o_N-o$ operations
of the above type as follows: We start with $(i,j)$ which is minimal w.r.t.\ "$\preceq$"
(because the corresponding process will increase the weight the least), perform
the operation on $(i,j)$ as many times as possible, say $s_{(i,j)}$ times,
then pass to the next pair $(i^\p,j^\p)\in I$ w.r.t.\ the order "$\preceq$" and so on,
until we have performed $s$ such processes in total.
Then, we arrive at a tuple $\ul{\rho}^\p$, such that the generic vanishing order of sections
from $W_{\ul{\rho}^\p}$ is precisely $o$. The difficult part is to show that
the corresponding weight will be, in fact, minimal.
\par
Fix a vanishing order $o<o_N$, let $w^{\ul{\alpha}}_{\iota,\ul{\rho}}(l)$
be the minimal weight of a section which vanishes of order $o$, and assume that $\ul{\rho}$ satisfies the inequalities of
Lemma~\ref{weightincrease1}. Then, we define natural numbers $s_{(i,j)}$ for
$(i,j)\in I$ inductively w.r.t.\ "$\preceq$" as follows:
For $(i,j)\in I$, set
\begin{eqnarray*}
c_i &=& \sum_{(i,j^\p)\prec (i,j)} s_{(i,j^\p)}
\\
c_j &=& \sum_{(i^\p,j)\prec (i,j)} s_{(i^\p,j)},
\end{eqnarray*}
where empty sums are by definition zero.
Then,
$$
s_{(i,j)}\q:=\q \min\bigl\{\,r_i-r_{i-1}-\rho_i-c_i,\rho_j-r_j+r_{j-1}-c_j\,\bigr\}.
$$
\begin{Obs}
\label{admissible}
{\rm i)} For every index $\kappa\in K^*$ and every admissible tuple $\ul{\rho}$,
satisfying the conditions of Lemma~\ref{weightincrease1}, we have
$$
\sum_{i\in K: i<\kappa}\bigl(\rho_i-r_i+r_{i-1}\bigr)
\q\le\q \sum_{j\in K^*: j\le \kappa} \bigl(r_j-r_{j-1}-\rho_j\bigr).
$$
From this, one easily infers that
$$
\begin{array}{cccl}
\rho_j &=& r_j-r_{j-1}+\sum_{(i,j)\in I} s_{(i,j)}&\hbox{for $j\in K$}
\\
\rho_i &=& r_i-r_{i-1}-\sum_{(i,j)\in I} s_{(i,j)}&\hbox{for $i\in K^*$},
\end{array}
$$
so that the $s_{(i,j)}$ determine $\ul{\rho}$.
\par
{\rm ii)} Suppose we are given a tuple $\ul{s}=(s^\p_{(i,j)},(i,j)\in I)$
with $0\le s_{(i,j)}^\p\le s_{(i,j)}$ for all $(i,j)\in I$.
Define $\ul{\rho}^{\ul{s}}$ by
$$
\begin{array}{cccl}
\rho^{\ul{s}}_j &=& r_j-r_{j-1}+\sum_{(i,j)\in I} s^\p_{(i,j)}&\hbox{for $j\in K$}
\\
\rho^{\ul{s}}_i &=& r_i-r_{i-1}-\sum_{(i,j)\in I} s^\p_{(i,j)}&\hbox{for $i\in K^*$}.
\end{array}
$$
The tuple $\ul{\rho}^{\ul{s}}$ is clearly admissible.
\end{Obs}
\begin{Lem}
\label{weightincrease2}
Fix a vanishing order $o<o_N$, and let $w^{\ul{\alpha}}_{\iota,\ul{\rho}}(l)$
be the minimal weight of a
section with vanishing order $o$ where $\ul{\rho}$ fulfills the conditions
of Lemma~{\rm\ref{weightincrease1}}. Then, $\ul{\rho}$ may be chosen in such a way that
the $s_{(i,j)}$ satisfy
$$
s_{(i,j)}\ =\ \min\left\{\,s-\sum_{(i^\p,j^\p)\prec(i,j)}s_{(i^\p,j^\p)},
 r_i-r_{i-1}-a_i+a_{i-1}-c_i, a_j-a_{j-1}-r_j+r_{j-1}-c_j\,\right\}.
$$
Here, $s:=\sum s_{(i,j)}$, and $c_i$ and $c_j$ are as before.
\end{Lem}
\begin{proof}
Assume that the assertion were wrong for $s_{(i,j)}$, i.e., $s_{(i,j)}$
is strictly smaller than the right hand side. Then, there are three cases.
In the first case, $s_{(i,j)}=r_i-r_{i-1}-\rho_i-c_i=\rho_j-r_j+r_{j-1}-c_j$.
In that case
\begin{equation}
\label{admissy}
\rho_i\ >\ a_i-a_{i-1}\q\hbox{and}\q \rho_j< a_j-a_{j-1}.
\end{equation}
Then, there is an $(i^\p,j^\p)\succ (i,j)$ with
$s_{(i^\p,j^\p)}>0$. Define
$\ul{\rho}^\p=(\rho^\p_1,...,\rho^\p_{k+1})$ with
$\rho^\p_{i^\p}:=\rho_{i^\p}+1$, $\rho^\p_{j^\p}:=\rho_{j^\p}-1$,
and $\rho^\p_j= \rho_j$ for $j\neq i^\p,j^\p$. The tuple
$\ul{\rho}^\p$ is still admissible, by \ref{admissible}. It is,
in fact, defined w.r.t.\ $\ul{s}$ with $s^\p_{(i,j)}=s_{(i,j)}$,
$(i,j)\neq (i^\p,j^\p)$, and
$s_{(i^\p,j^\p)}^\p=s_{(i^\p,j^\p)}-1$. Introduce
$\ul{\rho}^{\p\p}=(\rho^{\p\p}_1,...,\rho^{\p\p}_{k+1})$ by
$\rho^{\p\p}_i:=\rho^\p_i-1$, $\rho^{\p\p}_j:=\rho^\p_j+1$, and
$\rho^{\p\p}_{j^\p} :=\rho^\p_{j^\p}$ for $j^\p\neq i,j$. This is
again admissible, $w^{\ul{\alpha}}_{\iota,\ul{\rho}^{\p\p}}(l) \le
w^{\ul{\alpha}}_{\iota,\ul{\rho}}(l)$, and $W_{\ul{\rho}^{\p\p}}$
still contains sections of vanishing order $o$, by (\ref{admissy}). In
other words, we set $s^\p_{(i,j)}:=s_{(i,j)}+1$.
\par
In the second case, $s_{(i,j)}=r_i-r_{i-1}-\rho_i-c_i<\rho_j-r_j+r_{j-1}-c_j$.
Then, as $\rho_j-r_j+r_{j-1}-c_j>0$, there is an index $(i^\p,j)\succ(i,j)$
with $s_{(i^\p,j)}>0$. One may now proceed as before.
The last case, $s_{(i,j)}=\rho_j-r_j+r_{j-1}-c_j<r_i-r_{i-1}-\rho_i-c_i$, is handled
the same way.
\end{proof}
Given $s$, the condition in Lemma~\ref{weightincrease2} uniquely determines a
tuple $\ul{\rho}$ with $\sum s_{(i,j)}=s$ for which $w^{\ul{\alpha}}_{\iota,\ul{\rho}}(l)$
becomes minimal.
Note that $W_{\ul{\rho}}$ yields sections with vanishing order $\ge o_N-s$ where "$=$"
is achieved.
An immediate consequence is
\begin{Cor}
\label{weightincrease3}
{\rm i)}
Fix a vanishing order $o<o_N$, and let $w^{\ul{\alpha}}_{\iota,\ul{\rho}}(l)$ be the minimal
weight of a
section with vanishing order $o$ where $\ul{\rho}$ fulfills the conditions
of Lemma~{\rm\ref{weightincrease1}} and {\rm\ref{weightincrease2}}. Then,
$\sum s_{(i,j)}=o_N-o$.
\par
{\rm ii)} Denote by $w^{N,o}_{\iota,\ul{\alpha}}(l)$ the minimal weight of a
section with vanishing order $o$ at $N$. Then,
$$
w_{\iota,\ul{\alpha}}^{N,{o-1}}(l)-w_{\iota,\ul{\alpha}}^{N,o}(l) \q\ge\q
w_{\iota,\ul{\alpha}}^{N,{o}}(l)-w_{\iota,\ul{\alpha}}^{N,{o+1}}(l).
$$
\end{Cor}
Set ${\Bbb O}_o:=W_{\ul{\rho}}$,
where $\ul{\rho}$ is determined by the conditions of Lemma~\ref{weightincrease1} and
\ref{weightincrease2} and $\sum s_{(i,j)}=o_N-o$, and
${\Bbb O}:=\bigoplus_{o=0}^{o_N} {\Bbb O}_o$.
We have to find the minimal weights of sections in $H^0(C_\iota,\L_m^\iota)$ vanishing
of order $0\le o^\p\le m\cdot o_N-1$. We clearly have to look only at sections
in
$$
S^m {\Bbb O}=\bigoplus_{m_0,...,m_{o_N}:\sum m_\nu =m} S^{m_0}{\Bbb O}_0
\otimes\cdots\otimes S^{m_{o_N}} {\Bbb O}_{o_N}.
$$
Now, the sections in $S^{m_0}{\Bbb O}_0
\otimes\cdots\otimes S^{m_{o_N}} {\Bbb O}_{o_N}$
vanish of order at least
$m_1+2\cdot m_2+\cdots+o_N\cdot m_{o_N}$, and we can find some with exactly that
vanishing order. On the other hand, the weight of sections in that space
is
\begin{eqnarray*}
&&m_0\cdot w^{N,0}_{\iota,\ul{\alpha}}(l)+\cdots+m_{o_N}\cdot w^{N,o_N}_{\iota,\ul{\alpha}}(l)
\\
&=&m\cdot w^{N,0}_{\iota,\ul{\alpha}}(l)-
m_0\bigl(w^{N,0}_{\iota,\ul{\alpha}}(l)-w^{N,1}_{\iota,\ul{\alpha}}(l)\bigr)-\cdots-
m_{o_N}\bigl(w^{N,0}_{\iota,\ul{\alpha}}(l)-w^{N,o_N}_{\iota,\ul{\alpha}}(l)\bigr)
\\
&=& m\cdot w^{N,0}_{\iota,\ul{\alpha}}(l)-
m_0\bigl(w^{N,0}_{\iota,\ul{\alpha}}(l)-w^{N,1}_{\iota,\ul{\alpha}}(l)\bigr)-\cdots-
\\
&&\phantom{m\cdot w^{N,0}_{\iota,\ul{\alpha}}(l)}
- m_{o_N}\bigl((w^{N,0}_{\iota,\ul{\alpha}}(l)-w^{N,1}_{\iota,\ul{\alpha}}(l))
+\cdots+ (w^{N,o_N-1}_{\iota,\ul{\alpha}}(l)-w^{N,o_N}_{\iota,\ul{\alpha}}(l))\bigr).
\end{eqnarray*}
It follows easily from Corollary~\ref{weightincrease3} ii) that the
elements in $S^m{\Bbb O}$ producing sections of minimal weight vanishing of order  $o$ with
$(t-1)\cdot m\le o\le t\cdot m-1$ lie in
$$
\bigoplus_{i=1}^m S^{i}{\Bbb O}_{t-1}\otimes S^{m-i}{\Bbb O}_t,\qquad t=1,...,o_N.
$$
These contribute the weight
$$
m\cdot w^{N,t}_{\iota,\ul{\alpha}}(l)+\frac{m(m-1)}{2}\bigl(w^{N,t-1}_{\iota,\ul{\alpha}}(l)-
w^{N,t}_{\iota,\ul{\alpha}}(l)\bigr)
$$
to a basis for ${\Bbb H}:=H^0(C_{\iota},\L_m^\iota)/H^0\Bigl(\bigl(
\det(\E(l))_{|C_\iota}(-\sum_{s=1}^{\nu_{\iota}} o_{N_s} N_s)\bigr)^{\otimes m}\Bigr)$.
The total contribution to a basis for ${\Bbb H}$, coming from the node $N$, thus amounts to
$$
C^{N, \ul{\alpha}}_\iota(l,m)
\q:=\q m\bigl(w^{N,0}_{\iota,\ul{\alpha}}(l)+\cdots+w^{N,o_N}_{\iota,\ul{\alpha}}(l)\bigr)
+\frac{m(m-1)}{2}\bigl(w^{N,0}_{\iota,\ul{\alpha}}(l)-
w^{N,o_N}_{\iota,\ul{\alpha}}(l)\bigr).
$$
All in all, a basis for $H^0(C_{\iota},\L_m^\iota)$ will have minimal weight
$$
{C}^{\ul{\alpha}}_{\iota}(l,m)
\q:=\q
K^{\ul{\alpha}}_{\iota}(l,m)+\sum_{n=1}^{\nu_\iota} C^{N_n, \ul{\alpha}}_\iota(l,m).
$$
Let $N$ be an intersection of two components $C_{\iota}$ and $C_{\iota^\p}$ of $C$.
Then,
$$
w^{N,0}_{\iota,\ul{\alpha}}(l)\q=\q w^{N,0}_{\iota^\p,\ul{\alpha}}(l)\q=:\q
w^N_{\ul{\alpha}}(l).
$$
Let $\n$ be the set of nodes of $C$.
For large $m$, there is the exact sequence
$$
\begin{CD}
0 @>>> H^0(C,\det\E^{\otimes m}) @>>>\bigoplus_{\iota=1}^c
H^0(C_{\iota},\L_m^\iota) @>>> \bigoplus_{N\in\n} \C\cdot e_N@>>> 0.
\end{CD}
$$
This shows that $H^0(C,\det\E^{\otimes m})
=\bigoplus_{\iota=1}^c
H^0(C_{\iota},\L_m^\iota(-\sum_{n=1}^{\nu_\iota} N_n^{\iota}))
\oplus \bigoplus_{N\in\n} \C\cdot e_N$.
Thus, we see that the minimal weight of a basis for $H^0(C,\det \E^{\otimes m})$
is
$$
P^{\ul{\alpha}}_{\F_\bullet}(l,m)
\q:=\q \sum_{\iota=1}^c {C}^{\ul{\alpha}}_{\iota}(l,m)-m\cdot
\sum_{N\in\n}w^N_{\ul{\alpha}}(l).
$$
Note that this polynomial is intrinsically defined in terms
of the curve $C$, the filtration $\F_\bullet$, and $\ul{\alpha}$.
\subsubsection*{Case {\rm B)}: $C$ has nodal irreducible components}
In this case, we pass to the semistable curve $\pi\colon C^\p\lra C$,
where we introduce a projective line for every node at which $C$ is irreducible
and the filtration by the $\F_j$ is not a filtration by subbundles,
and pull-back $\E$ to $\E^\p$ on $C^\p$.
Now, let $C^\p_{\iota}$ be any irreducible component of $C^\p$.
Then, $\nu_\iota:=\pi_{|C_\iota^\p}\colon C^\p_\iota\lra C_\iota$ is a partial normalization.
For any node $N\in C_\iota$ which is resolved by $\nu_\iota$, $W^\iota_{\min}$ will
produce sections of $\det(\E^\p(l))$ which vanish at both points $p_{N,1}$
and $p_{N,2}$ in $\nu^{-1}_\iota(N)$.
As the space of sections of $\det(\E(l))$ vanishing at $N$ identifies with the space
of sections of $\det(\E^\p(l))$ vanishing at both $p_{N,1}$
and $p_{N,2}$, it is easy to see that
the analogs of Lemma~\ref{generate} and Corollary~\ref{asi1} continue to hold.
The rest of the considerations clearly go through as before.
\subsubsection*{H-semistability}
We are now ready to define our semistability concept.
\begin{Def}
\label{H-stab}
A pair $(C,\E)$, consisting of a semistable curve $C$ and a vector bundle $\E$ of rank $r$
on $C$ with $\chi(\E)=\chi$,
will be called \it H-(semi)stable\rm, if it satisfies the following conditions
\begin{enumerate}
\item The push forward $\pi_*(\E)$ to the stable model via $\pi\colon C\lra
      C^\p$ is a semistable torsion free sheaf (see the remark in the introduction).
\item For every maximal filtration of $\F_\bullet$ of $\pi_*(\E)$ by destabilizing
      subsheaves and every vector $\ul{\alpha}$ of non-negative rational numbers, there
      exists an index $l^*$, such that for
      all $l\ge l^*$
      $$
      P^{\ul{\alpha}}_{\F_\bullet}(l,m)\ (\preceq)\ 0\q\hbox{as polynomial in $m$}.
      $$
      \end{enumerate}
\end{Def}
Note that this semistability concept has all the properties that were asserted in the
introduction.
Part ii) of the Main Theorem is a direct consequence of
\begin{Thm}
There exist an index $l_0$ and for every $l\ge l_0$ an index $m(l)$, such that for
every $l\ge l_0$, $m^\p\ge m(l)$, and
every pair $(C,\E)$, consisting of a semistable curve $C$ and a vector bundle $\E$ of rank $r$
on $C$ with $\chi(\E)=\chi$, which satisfies {\rm (1)} and {\rm (2)} of {\rm \ref{H-stab}}
$$
P^{\ul{\alpha}}_{\F_\bullet}(l,m^\p)\ (\le)\ 0
\q\Longleftrightarrow\q
P^{\ul{\alpha}}_{\F_\bullet}(l,m)\ (\preceq)\ 0\q\hbox{as polynomial in $m$}
$$
for every filtration $\F_\bullet$ and every tuple $\ul{\alpha}$.
\end{Thm}
\begin{proof}
First note that, given $\ul{\alpha}$, $P^{\ul{\alpha}}_{\F^\bullet}$ depends only on the following data
\begin{itemize}
\item The tuples $(\rk \F^\iota_1,...,\rk\F^\iota_k)$, $\iota=1,...,c$. These determine
      all the Hilbert polynomials of the $\F^\p_j$, $j=1,...k$, because these are
      destabilizing sheaves.
\item The tuples $(a^N_1,...,a^N_k)$, $N$ a node of $C$, and $a^N_j$ the dimension
      of the image of $\F_j$ in the fibre of $\E$ at $N$, $j=1,...,k$.
\end{itemize}
By boundedness, the sets of data of the above type is in fact
finite. Therefore, we will be done, once we have shown that, for
a given set of such data, we have to take only finitely many
vectors $\ul{\alpha}$ into account.
\par
Given tuples $(r_1^\iota,...,r_k^\iota)$, $\iota=1,...,c$, and $(a^N_1,...,a^N_k)$, $N$
a node of $C$, we define sets $K_{N,\iota}$ and $K^*_{N,\iota}$ as before.
Note that in our construction before, we had to look at the quantities
$\gamma_j(l)-\gamma_i(l)$, $(i,j)\in K^*_{N,\iota}\times K_{N,\iota}$.
By Formula~(\ref{weight}),
$$
\gamma_j(l)-\gamma_i(l)\q=\q \sum_{t=i}^{j-1} \alpha_t\cdot\chi_l.
$$
For every ordering "$\preceq_{N,\iota}$" of $K^*_{N,\iota}\times K_{N,\iota}$, we
get the set of inequalities
$$
(*)_{\preceq_{N,\iota}}\colon\q\q
\sum_{t=i^\p}^{j^\p-1} \alpha_t\q \le\q \sum_{t=i}^{j-1} \alpha_t,\qquad
(i^\p,j^\p)\preceq_{N,\iota} (i,j).
$$
Let $Q\subset \R^k$ be the quadrant of vectors all the entries of which are non-negative. This is
a rational polyhedral cone. For a given ordering "$\preceq_{N,\iota}$", the
inequalities $(*)_{\preceq_{N,\iota}}$ define a proper rational polyhedral subcone of
$Q$. Given two distinct orderings, the resulting cones will meet only along faces, i.e.,
if we let "$\preceq_{N,\iota}$" vary over all possible orderings, we get a fan decomposition
$Q=\bigcup_{\beta=1}^{B_{N,\iota}} Q^{N,\iota}_\beta$ of $Q$.
\par
We have seen that once the ordering "$\preceq_{N,\iota}$" is fixed,
every given vanishing  order $o$ uniquely determines a vector $\ul{\rho}$
with $w^{\ul{\alpha}}_{\iota,\ul{\rho}}(l)=w^{N,o}_{\iota,\ul{\alpha}}(l)$ for all
$\ul{\alpha}$ in the cone $Q^{N,\iota}_\beta$ cut out by the inequalities
$(*)_{\preceq_{N,\iota}}$.
In particular, for $\ul{\alpha},\ul{\alpha}^\p\in Q^{N,\iota}_\beta$
\begin{equation}
\label{add1}
w^{N,o}_{\iota,\ul{\alpha}+\ul{\alpha}^\p}(l)
\q=\q
w^{N,o}_{\iota,\ul{\alpha}}(l)+w^{N,o}_{\iota,\ul{\alpha}^\p}(l)
\end{equation}
for all possible vanishing orders.
\par
As the intersection of two rational polyhedral cones is again a rational polyhedral
cone, we can form the rational polyhedral cones of the form
$Q^{N_1,\iota_1}_{\beta_1}\cap\cdots\cap Q^{N_\nu,\iota_\nu}_{\beta_\nu}$.
Here, $N_1,...,N_\nu$ are the nodes of $C$
(or, in Case B), the nodes of the corresponding partial normalization),
$\iota_i$ is an index such that $C_{\iota_i}$
contains $N_i$, and $\beta_i\in\{\,1,..., B_{N,\iota_i}\,\}$, $i=1,...,\nu$.
This defines a fan decomposition $Q=\bigcup_{\beta=1}^B Q_B$.
\par
Let $\F_\bullet$ be a maximal filtration, realizing the data
$(r_1^\iota,...,r_k^\iota)$, $\iota=1,...,c$, and $(a^N_1,...,a^N_k)$, $N$
a node of $Q$. Then, for every $\beta\in\{\, 1,...,N\,\}$, and any
$\ul{\alpha},\ul{\alpha}^\p\in Q_\beta$, we have, by (\ref{add1}),
\begin{equation}
\label{add2}
P_{\F^\bullet}^{\ul{\alpha}+\ul{\alpha}^\p}(l,m)
\q=\q
P_{\F^\bullet}^{\ul{\alpha}}(l,m)
+
P_{\F^\bullet}^{\ul{\alpha}^\p}(l,m).
\end{equation}
For every edge $e$ of the cone $Q_\beta$, denote the minimal integral generator by
$\ul{\alpha}_{e,\beta}$. Then, by (\ref{add2}), we have to verify the inequalities
in the Definition~\ref{H-stab} only for $\ul{\alpha}$ in the finite
set $\{\,\ul{\alpha}_{e,\beta}\,|\,\beta=1,...,B,\ e \hbox{ an edge of $Q_\beta$}\,\}$.
The theorem is now settled.
\end{proof}
\subsection{The Hilbert compactification as a moduli space}
Introduce the functors
$$
\ul{\rm HC}^{\rm (s)s}(g;\chi,r)\colon \ul{\rm Schemes}_\C \lra
\ul{\rm Sets}
$$
which assign to every scheme $S$ the equivalence classes of pairs
$({\cal C}_S, \E_S)$ where $\pi\colon{\cal C}_S\lra S$ is a flat
family of semistable curves, and ${\E}_S$ is an $S$-flat sheaf,
such that, for every closed point $s\in S$, the restriction
$\E_{S|\pi^{-1}(s)}$ is an H-(semi)stable vector bundle of
uniform rank $r$ and Euler characteristic $\chi$. Two families
$({\cal C}_S, \E_S)$ and $({\cal C}^\p_S, \E^\p_S)$ are \it
equivalent\rm, if there are an isomorphism $\phi_S\colon {\cal
C}_S^\p\lra {\cal C}_S$ and a line bundle $\L_S$ on $S$, such that
$$
\phi_S^*\bigl(\E_S\otimes \pi^{*}(\L_S)\bigr)\q\cong\q \E_S^\p.
$$
Then, our considerations imply
\begin{Thm}
{\rm i)} There is a natural transformation $\theta\colon \ul{\rm
HC}^{\rm ss}(g;\chi,r) \lra h_{\frak H\frak C(g;\chi,r)}$, such
that for every other scheme ${\frak S}$ and every other natural
transformation $\theta^\p\colon \ul{\rm HC}^{\rm
ss}(g;\chi,r)\lra h_{\frak S}$, one has a unique morphism
$t\colon {\frak H\frak C}(g;\chi,r)\lra {\frak S}$ with
$\theta^\p=h(t)\circ \theta$.
\par
{\rm ii)} The space ${\frak H\frak C(g;\chi,r)}$
contains an open subscheme ${\frak H\frak C(g;\chi,r)}^{\rm s}$ which is a coarse moduli scheme
for $\ul{\rm HC}^{\rm s}(g;\chi,r)$.
\end{Thm}
\section{Properties of the Hilbert compactification}
\subsection{Dimension and smooth points}
\label{smoothsection}
We will call a pair $(C,\E)$ with $C$ a semistable curve and $\E$
a vector bundle on $C$ \it strictly H-stable\rm, if it is
H-stable and there is no automorphism $\phi\colon C\lra C$ with
$\phi^*\E\cong\E$.
\begin{Rem}
\label{Simple} If $(C,\E)$ is strictly H-stable, then $\E$ must be
a simple bundle, i.e., $\End(\E)\cong \C\cdot {\rm id}_\E$. In
fact, the universal bundle on $\widehat{\frak H}:=\widehat{\frak
H}(g;\chi_l,r)$ possesses a $\GL(V^{\chi_l})$-linearization whence
the $\GL(V^{\chi_l})$-stabilizers of a point in $\widehat{\frak
H}$ corresponding to a strictly H-stable pair $(C,\E)$ identify
with the automorphisms of $\E$ on $C$ which form a dense set in
the space of endomorphisms. Therefore, $\End(\E)$ can have
dimension at most one, because the $\GL(V^{\chi_l})$-stabilizer
may have dimension at most one.
\end{Rem}
Let ${\frak H\frak C}(g;\chi,r)^\star\subset {\frak H\frak
C}(g;\chi,r)$ be the open subset parameterizing the strictly
H-stable curves.
\begin{Thm}
\label{smooth} {\rm i)} The Hilbert compactification is a normal
variety of dimension $3g-3+r^2(g-1)+1$.
\par
{\rm ii)} The subset ${\frak H\frak C}(g;\chi,r)^\star$ is smooth.
\end{Thm}
\begin{proof}
Let $\kappa_g\colon {\frak H\frak C}(g;\chi,r)\lra \overline{\frak
M}_g$ be the natural morphism. The irreducibility and the
dimension statement in i) are clear, because they are known for
the preimage ${\cal U}$ of the moduli space ${\frak M}_g$ of
smooth curves under $\kappa_g$ (see \cite{Pa}), and we have seen
in \ref{Semstabred} that ${\cal U}$ is dense in the Hilbert
compactification.
\par
For the remaining statements, let ${\frak H}^0(g;\chi,r)$ be as in
the introduction, and ${\frak H}^\dagger(g;\chi,r)\subset {\frak
H}^0(g;\chi,r)$ the open subset of those $(C,q_C\colon
V^\chi\otimes{\Oh}_C\lra \E)$ for which $H^1(\E_C)$ vanishes. As our
considerations in Chapter 2 show, the Hilbert compactification is
a quotient of an open subset of ${\frak H}^\dagger(g;\chi,r)$.
Please accept for the moment the following statement.
\begin{Prop}
\label{really smooth} The scheme ${\frak H}^\dagger(g;\chi,r)$ is
smooth.
\end{Prop}
This proposition settles i). Let ${\frak H}^\star(g;\chi,r)$ be
the open part of the Hilbert scheme which parameterizes the
strictly H-stable objects. Now, statement ii) follows, because the
quotient morphism ${\frak H}^\star(g;\chi,r)\lra {\frak H\frak
C}(g;\chi,r)^\star$ is a principal $\PGL(V^{\chi_l})$-bundle.
\par
We now turn to the proof of Proposition~\ref{really smooth}.
Let ${\frak M}_g^\star$ be the moduli space of automorphism free
smooth curves, and set
$$
{\cal U}^\star\q:=\q \kappa_g^{-1}\bigl({\frak
M}_g^\star\bigr)\cap {\frak H\frak C}^{\rm s}.
$$
Then, ${\cal U}^\star$ is a smooth quasi-projective variety of
dimension $3g-3+r^2(g-1)+1$. Since the quotient morphism is over
${\cal U}^\star$ a principal $\PGL(V^{\chi_l})$-bundle, the
preimage of ${\cal U}^\star$ under the quotient morphism is a
smooth quasi-projective variety of dimension
$3g-3+r^2(g-1)+\chi_l^2$. Moreover, by~\ref{Semstabred}, it is
dense in ${\frak H}^\dagger(g;\chi,r)$, whence the latter is an
irreducible scheme of the same dimension. To prove smoothness, we
have to determine the dimension of the tangent spaces. If $x\in
{\frak H}^\dagger(g;\chi,r)$ corresponds to the curve
$C_x\hookrightarrow {\frak G}$, the tangent space to ${\frak
H}^\dagger(g;\chi,r)$ at $x$ is given by $\Hom\bigl({\cal
I}_{C_x}/{\cal I}_{C_x}^2,{\Oh}_{C_x}\bigr)$. Since $C$ is a local
complete intersection (which is an intrinsic property by \cite{LS}, Prop.\ 3.2.1 and Cor.\ 3.2.2), 
the conormal sheaf ${\cal I}_{C_x}/{\cal
I}_{C_x}^2$ is locally free, and we  have the exact sequence
$$
\begin{CD}
0 @>>> {\cal I}_{C_x}/{\cal I}_{C_x}^2 @>>> {\Oh}mega^1_{\frak G|C_x}
@>>> {\Oh}mega^1_{C_x} @>>> 0.
\end{CD}
$$
Here, the left exactness follows, because (a) the sequence is in
any case exact away from the nodes of $C_x$ and (b) since ${\cal
I}_{C_x}/{\cal I}_{C_x}^2$ is torsion free, it does not contain
any subsheaf the support of which has dimension strictly less
than one. We derive the exact sequence
\begin{eqnarray*}
0&\lra& \Hom({\Oh}mega_{C_x}^1,{\Oh}_{C_x})\lra H^0(T_{\frak
G|C_x})\q\lra
\\
& \lra & \Hom\bigl({\cal I}_{C_x}/{\cal
I}_{C_x}^2,{\Oh}_{C_x}\bigr)\lra {\rm
Ext}^1({\Oh}mega^1_{C_x},{\Oh}_{C_x})\q \lra\q H^1(T_{\frak G|C_x}).
\end{eqnarray*}
We claim that $H^1(T_{\frak G|C_x})$ vanishes. For this, we use
the exact sequence
$$
\begin{CD}
0 @>>> \E nd(\E_{\frak G}) @>>> \E_{\frak G}^{\oplus \chi_l} @>>>
T_{\frak G}@>>> 0.
\end{CD}
$$
Let $\E$ be the restriction of $\E_{{\frak G}}$ to the curve
$C_x$, so that we obtain the exact sequence
$$
\begin{CD}
0 @>>> \E nd(\E) @>>> \E^{\oplus \chi_l} @>>> T_{\frak G|C_x}@>>>
0.
\end{CD}
$$
Now, our assumption is that $H^1(\E)^{\oplus
\chi_l}=H^1(\E^{\oplus \chi_l})$ vanishes, and, for dimension
reasons, $H^2(\E nd(\E))=0$, whence also $H^1(T_{{\frak
G}|C_x})=0$, as asserted. We also see that
$$
h^0(T_{\frak G|C_x})\q=\q \chi_l\cdot h^0(\E) -\chi(\E nd(\E))
\q=\q\chi_l^2+r^2(g-1).
$$
Next, by Serre duality
$$
\dim\bigl(\Hom({\Oh}mega^1_{C_x},{\Oh}_{C_x})\bigr) - \dim\bigl({\rm
Ext}^1({\Oh}mega_{C_x}^1,{\Oh}_{C_x})\bigr)=\chi({\Oh}mega_{C_x}^1\otimes
\omega_{C_x})=3g-3.
$$
The exact sequence above thus shows
$$
\dim\bigl(\Hom({\cal I}_{C_x}/{\cal I}_{C_x}^2,{\Oh}_{C_x})\bigr) \
=\ h^0(T_{\frak G|C_x})+\chi({\Oh}mega_{C_x}^1\otimes \omega_{C_x})\
=\  \chi_l^2+r^2(g-1)+3g-3.
$$
This proves that $x$ is a smooth point of ${\frak
H}^\dagger(g;\chi,r)$.
\end{proof}
\begin{Rem}
\label{deformations}
i) Deligne and Mumford \cite{DM} applied
Schlessinger's deformation theory \cite{Schl} in order to show
that any semistable curve $C$ admits a \it miniversal deformation
\rm over the base scheme ${\cal M}:={\rm Spec} \C\formall
t_1,...,t_N\formalr$ with $N:={\rm dim}\bigl({\rm
Ext}^1({\Oh}mega_C,{\Oh}_C)\bigr)$. This means that there is a family
${\cal C}_{\cal M}\lra {\cal M}$ of curves parameterized by ${\cal
M}$ with $C$ as the fibre over the origin, such that for any flat
family of curves ${\cal C}_B\lra B$ with $B$ the spectrum of a
local Artin algebra and $C$ as the fibre over the closed point
there is a morphism $\phi\colon B\lra {\cal M}$ with ${\cal
C}_B\cong{\cal C}_{\cal M}\times_{\cal M} B$. Moreover, $\phi$ is
unique in case $B={\rm
Spec}\bigl(\C[\eps]/\langle\eps^2\rangle\bigr)$. The tangent space
to $0\in {\cal M}$ thus identifies with ${\rm
Ext}^1({\Oh}mega_C,{\Oh}_C)$. Finally, suppose $C$ has $M$ nodes, then
the space of \it local deformations \rm is the deformation space
of this set of nodes and, thus, identifies with ${\cal M}_{\rm
loc}={\rm Spec}\C\formall u_1,...,u_M\formalr$. Here, one can
arrange the generators $t_1,...,t_N$ and $u_1,...,u_M$ in such a
way that the natural morphism ${\cal M}\lra {\cal M}_{\rm loc}$
comes from the homomorphism $\C\formall
 u_1,...,u_M\formalr \lra
\C\formall t_1,...,t_N\formalr$, $u_i\lma t_i$, $i=1,...,M$.
\par
Now, let $x\in {\frak H}^\dagger(g;\chi,r)$, and let ${\frak U}$
be its formal neighborhood. By the smoothness of ${\cal M}$ and
its versality, the universal curve over the Hilbert scheme
${\frak H}^\dagger(g;\chi,r)$ provides us with a morphism
$\phi\colon {\frak U}\lra {\cal M}$ the differential of which is
the map
$$
\Hom\bigl({\cal I}_{C_x}/{\cal I}_{C_x}^2,{\Oh}_{C_x}\bigr)\lra {\rm
Ext}^1({\Oh}mega^1_{C_x},{\Oh}_{C_x}).
$$
As we have seen before, this map is surjective, so that $\phi$ is
a submersion whence a smooth morphism.
\par
ii) If we work in the setting of Artin and Deligne-Mumford stacks,
we can sharpen the second statement of Theorem~\ref{smooth}. We
will do this in Section~\ref{Stack} below.
\end{Rem}
\subsection{Existence of universal families}
\label{UnivFam}
The aim of this section is to prove
\begin{Thm}
\label{Poin}
Suppose $\chi$ and $r$ are coprime. Then, every point
$x\in {\frak H\frak C}(g;\chi,r)^\star$ possesses an \'etale
neighborhood $U$, such that there exists a universal family over
$U$.
\end{Thm}
\subsubsection{The universal curve over ${\frak H\frak C(g;\chi,r)}^{\rm s}$}
Let ${\cal C}_{\widehat{\frak H}}\hookrightarrow {\widehat{\frak
H}}\times {\frak G}$ be the universal closed subscheme. It is
clearly invariant under the $\SL(V^{\chi_l})$-action on
${\widehat{\frak H}}\times {\frak G}$. On ${\widehat{\frak H}}$,
we choose a line bundle ${\Oh}_{\widehat{\frak H}}(1):={\frak
L}_{a,m}$ where $m$ and $a\gg 0$ are chosen in such a way that
the conclusion of the Main Theorem holds, and on ${\frak G}$ the
usual ample line bundle ${\Oh}_{\frak G}(1)$. Then, for positive
integers $s$ and $t$,
$$
{\cal L}_{s,t}:=\Bigl(\pi_{\widehat{\frak
H}}^*{\Oh}_{{\widehat{\frak H}}}(s)\otimes \pi_{\frak G}^*{\Oh}_{\frak
G}(t)\Bigr)_{|{\cal C}_{\widehat{\frak H}}}
$$
is an ample $\SL(V^{\chi_l})$-linearized line bundle on ${\cal
C}_{\widehat{\frak H}}$. If we choose $s/t$ large enough,
Proposition~\ref{suitable} grants that points in the preimage
${\cal C}^{\rm s}_{\widehat{\frak H}}$ of the points in
${\widehat{\frak H}}$ which are stable w.r.t.\ the linearization
in ${\frak L}_{a,m}$ are stable w.r.t.\ the linearization in
${\cal L}_{s,t}$. Therefore, the geometric quotient
$$
{\cal C}_{{\frak H\frak C}^{\rm s}} \q:=\q {\cal
C}_{\widehat{\frak H}}^{\rm s}\catqot \SL(V^{\chi_l})
$$
exists. Here, we have set ${\frak H\frak C}^{\rm s}:={\frak H\frak
C(g;\chi,r)}^{\rm s}$. Moreover, there is a natural morphism
$$
\sigma\colon\  {\cal C}_{{\frak H\frak C}^{\rm s}}\lra {\frak
H\frak C}^s.
$$
We call ${\cal C}_{\frak H\frak C^s}$ --- abusively
--- the \it universal curve\rm.
For an H-stable pair $(C,\E)$, define
$$
{\rm Aut}(C,\E)
\q:=\q\bigl\{\, \alpha\colon C\stackrel{\cong}{\lra} C\,|\,
\alpha^*\E\cong\E\,\bigr\}.
$$
From the GIT set up, it follows that the
$\PGL(V^{\chi_l})$-stabilizer of a stable point $(C,q\colon
V^{\chi_l}\otimes{\Oh}_C\lra\E)$ in ${\widehat{\frak H}}$ identifies
with the group ${\rm Aut}(C,\E)$. Thus, we see
\begin{Cor}
For any point $[C,\E]\in {\frak H\frak C}^s$, the fibre
$\sigma^{-1}[C,\E]$ is isomorphic to the curve $C/{\rm
Aut}(C,\E)$.
\end{Cor}
In particular, we may hope for a universal family only over the
open subset ${\frak H\frak C}(g;\chi,r)^\star$.
\subsubsection{Proof of Theorem~\ref{Poin}}
Now, let $x\in {\frak H\frak C}^\star:={\frak H\frak
C}(g;\chi,r)^\star$. Then, $x$ has an \'etale neighborhood $U$,
such that the family $\sigma_U\colon {\cal C}_U:={\cal C}_{\frak
H\frak C^\star}\times_{\frak H\frak C^\star} U\lra U$ possesses a
section which meets every fibre in a smooth point. Then,
$\L_{U}:={\Oh}_{{\cal C}_U}\bigl(\sigma_U(U)\bigr)$ is a relative
ample invertible sheaf. Let ${\frak U}^\star\subset
{\widehat{\frak H}}$ be the $\SL(V^{\chi_l})$-invariant open
subset which parameterizes the strictly H-stable points, and let
$\psi\colon {\frak U}^\star\lra {\frak H\frak C}^\star$ be the
quotient morphism. Then, since Mumford's GIT supplies universal
geometric quotients, the map $\psi_U\colon {\cal U} :={\frak
U}^\star\times_{\frak H\frak C^\star} U\lra U$ is a geometric
quotient, too. For the same reason, the vertical maps in the
following Cartesian diagram are both geometric quotients:
$$
\begin{CD}
{\cal C}_{\cal U}:={\cal C}_{\frak H\frak C^\star}\times_{\frak
H\frak C^\star} {\cal U} @>\sigma_{\cal U} >> {\cal U}
\\
@V\psi_{{\cal C}_{ U}}VV @VV\psi_U V
\\
{\cal C}_{U} @>>> U.
\end{CD}
$$
Define $\L_{\cal U}$ as the pullback of $\L_U$ under the quotient
morphism $\psi_{{\cal C}_{U}}$. This is obviously an
$\SL(V^{\chi_l})$-linearized relative ample invertible sheaf.
\par
We have seen that the $\SL(V^{\chi_l})$-stabilizer of a point
$u\in {\cal U}$ corresponding to a strictly H-stable object
consists exactly of the scalar matrices (Remark~\ref{Simple}).
The same holds for the points $x\in {\cal C}_{\cal U}$. Let
$\mu_{\chi_l} \subset \C^*$ be the subgroup generated by a
primitive $\chi_l$-th root of unity $\zeta$. For any point $u\in
{\cal U}$ or $x\in {\cal C}_{\cal U}$, the
$\SL(V^{\chi_l})$-stabilizer now identifies with $\mu_{\chi_l}$.
If $\F$ is an $\SL(V^{\chi_l})$-linearized sheaf on ${\cal U}$ or
${\cal C}_{\cal U}$, we will say that it is of \it weight $k$\rm,
if $\mu_{\chi_l}$ acts by $\zeta^k\cdot {\rm id}_{\F\langle
x\rangle}$ for all $x\in {\cal U}$ or $x\in {\cal C}_{\cal U}$,
respectively. For example, $\L_{\cal U}$ is of weight zero.
\par
By construction, we have the quotient $q_{\cal U} \colon
V^{\chi_l}\otimes {\Oh}_{\cal U} \lra \E_{\cal U}$, and the question
we have to answer is whether $\E_{\cal U}$ descends --- possibly
after tensorizing it with the pullback of a line bundle on ${\cal
U}$ --- to the quotient ${\cal C}_U$. By Kempf's descent lemma
(see \cite{DN}), an $\SL(V^{\chi_l})$-linearized vector bundle
$\E$ on ${\cal C}_{\cal U}$ descends to the quotient, if and only
if it is of weight zero. Now, the sheaf $\E_{\cal U}$ is of
weight one. Thus, our task will be to find an
$\SL(V^{\chi_l})$-linearized invertible sheaf ${\cal A}_{\cal U}$
of weight one. Then, $\E_{\cal U}\otimes \sigma_{\cal
U}^*\bigl({\cal A}^\vee_{\cal U}\bigr)$ will descend to ${\cal
C}_U$, and we will be done.
\par
For any $m$, the sheaf $\E_{\cal U}[m]:= \E_{\cal U}\otimes
\L_{\cal U}^{\otimes m}$ is an $\SL(V^{\chi_l})$-linearized vector
bundle of weight one, and, if $m$ is sufficiently large,
$$
\F_m:=\sigma_{\cal U*} \bigl(\E_{\cal U}[m]\bigr)
$$
will be an $\SL(V^{\chi_l})$-linearized vector bundle of rank
$\chi_l+r\cdot m$ and weight one. Then, for $m\gg 0$,
$$
{\cal N}_{\cal U}\q:=\q \det(\F_{m+1})\otimes \det(\F_{m})^\vee
$$
is a line bundle of weight $r$. Since $\chi$ and $r$ and thus also
$c:=\chi_l+r\cdot m$ and $r$ are coprime, we may find integers
$\alpha$ and $\beta$ with $\alpha\cdot c+\beta\cdot r=1$, so that
$$
{\cal A}_{\cal U}\q:=\q \det(\F_m)^{\otimes \alpha}\otimes {\cal
N}_{\cal U}^{\otimes \beta}
$$
will indeed have weight one. \qed
\subsection{The fibres of the morphism  $\kappa_g$}
As before, let $\kappa_g\colon {\frak H\frak C}(g;\chi,r)\lra
\overline{\frak M}_g$ be the natural morphism. By $\ol{\frak
M}^\star_g$, we denote the quasi-projective moduli space of
automorphism free stable curves. In this section, we want to
establish
\begin{Thm}
\label{fibres} For any stable curve $C_0$ without automorphisms,
the variety
$$
\kappa_g^{-1}\bigl([C_0]\bigr)\cap {\frak H\frak C}^{\rm
s}(g;\chi,r)
$$
has only analytical normal crossings as singularities.
\end{Thm}
We will follow the strategy of Gieseker's paper \cite{Gie} in
order to prove the result.
\subsubsection{A family of semistable curves with fixed stable
model}
Let $C_0$ be a fixed stable curve. Then, by the results of
Deligne and Mumford \cite{DM}, $C_0$ has a universal deformation
over ${\cal M}:={\rm Spec} \C\formall t_1,...,t_N\formalr$ with
$N:={\rm dim}\bigl({\rm Ext}^1({\Oh}mega_{C_0},{\Oh}_{C_0})\bigr)=3g-3$.
Moreover, let ${\cal M}_{\rm loc}={\rm Spec}\C\formall
u_1,...,u_M\formalr$ be the deformation space of the nodes of
$C_0$. Finally, there is the morphism ${\cal M}\lra {\cal M}_{\rm
loc}$ normalized in such a way that $\C\formall
 u_1,...,u_M\formalr \lra
\C\formall t_1,...,t_N\formalr$ is given by $u_i\lma t_i$,
$i=1,...,M$, and $t_i=0$ is the equation of the $i$-th node of
$C_0$, $i=1,...,M$.
\par
Next, let $C$ be a semistable curve the stable model of which is
$C_0$, and let $\pi\colon C\lra C_0$ be the contraction map. Let
$R_1,...,R_S$ be the maximal connected chains of rational curves
which are contracted by $\pi$. We label the nodes $c_1,...,c_M$ of
$C_0$ in such a way that $\{c_i\}=\pi(R_i)$, $i=1,...,S$. We then
define
$$
\n\q:=\q \bigl\{\, t_i=0,\ i=S+1,...,N\,\bigr\}\q\subset\q {\cal
M}.
$$
By restriction of the universal family over $\M$, we find a family
$\sigma_\n\colon {\cal C}_{\n}\lra \n$ which is smooth outside the
nodes "which don't move", i.e., outside the nodes
$c_{S+1},...,c_M$. For $i=1,...,S$, let $d_{i,j}$,
$j=1,...,\iota_i$, be the nodes of $C$ mapping under $\pi$ to
$c_i$. Define
$$
{\cal Q}\q:=\q \Spec\bigl(\C\formall\, x_{i,j};\ i=1,...,S,
j=1,...,\iota_i\,\formalr\bigr).
$$
The homomorphism
\begin{eqnarray*}
\phi^*\colon\  \C\formall\, t_1,...,t_S\,\formalr &\lra &
\C\formall\,  x_{i,j};\ i=1,...,S, j=1,...,\iota_i\,\formalr
\\
t_i &\lma & x_{i,1}\cdot ...\cdot x_{i,\iota_i},\q i=1,...,S,
\end{eqnarray*}
defines a morphism $\phi\colon {\cal Q}\lra {\cal N}$. The pull
back of the family  ${\cal C}_\n$ provides us with the family
$\sigma_{\cal Q}\colon {\cal C}_{\cal Q}\lra {\cal Q}$ of stable
curves. Near the $i$-th node, the family ${\cal C}_{\cal Q}$ is
defined by the equation
$$
y_i\cdot z_i - x_{i,1}\cdot ...\cdot x_{i,\iota_i}\q =\q 0,
$$
for appropriate parameters $y_i$ and $z_i$, $i=1,...,S$.
\par
Now, let
$$
\theta\colon \q \widehat{\cal C}_{\cal Q} \lra {\cal Q}
$$
be the blow up of the curve ${\cal C}_{\cal Q}$ along the ideal
generated by $y_1$ and $x_{1,1}$. Near the node $c_1$, we may
embed ${\cal C}_{\cal Q}$ into
$$
{\Bbb A}\q=\q \Spec\bigl(\C\formall\, y_1,z_1,x_{i,j};\ i=1,...,S,
j=1,...,\iota_i\,\formalr\bigr).
$$
The blow up $\widehat{\Bbb A}$ of ${\Bbb A}$ along $y_1$ and
$x_{1,1}$ is the scheme
$$
\widehat{\Bbb A}\q:=\q \Bigl\{\, \bigl((y_1,z_1,x_{i,j};\
i=1,...,S), [w_0:w_1]\bigr)\, \big|\, y_1\cdot w_0 = x_{1,1}\cdot
w_1\,\Bigr\}\q\subset\q {\Bbb A}\times{\Pe}_1.
$$
One checks that the strict transform of ${\cal C}_{\cal Q}$ is
given in the chart $w_0=1$ by the equation
$$
w_1\cdot z_1- x_{i,2}\cdot ...\cdot x_{i,\iota_1}\q=\q 0
$$
and in the chart $w_1=1$ by
$$
w_0\cdot y_1- x_{i,1}\q=\q 0.
$$
We may now iterate the blow up, i.e., blow up $\widehat{\cal
C}_{\cal Q}$ at the ideal generated by $w_1$ and $x_{i,2}$ and so
on and perform the same procedure at the other nodes, too, in
order to construct a flat family
$$
\widetilde{\sigma}_{\cal Q}\colon\  \widetilde{\cal C}_{\cal Q}
\lra {\cal Q}
$$
with $C$ as the fibre over the origin. By construction,
$\widetilde{\cal C}_{\cal Q}$ is given near the node $d_{i,j}$ by
the equation
$$
y_{i,j}\cdot z_{i,j}- x_{i,j}\q=\q 0
$$
for suitable local parameters $y_{i,j}$ and $z_{i,j}$,
$i=1,...,S$, $j=1,...,\iota_i$. In particular, it is near
$d_{i,j}$ isomorphic to the miniversal deformation of that node,
and $x_{i,j}=0$ is the locus where the node $d_{i,j}$ "is kept",
$i=1,...,S$, $j=1,...,\iota_i$.
\par
By ${\cal X}\hookrightarrow{\cal Q}$, we denote the subscheme
defined by the equations
$$
x_{i,1}\cdot ...\cdot x_{i,\iota_i}=0,\q i=1,...,S.
$$
The scheme ${\cal X}$ obviously has only analytical normal
crossing singularities.
\subsubsection{The versality property of $\widetilde{\cal C}_{\cal
Q}$}
The family $\widetilde{\sigma}_{\cal Q}\colon \widetilde{\cal C
}_{\cal Q}\lra {\cal Q}$ together with the ${\cal Q}$-morphism
$\pi_{\cal Q}\colon \widetilde{\cal C}_{\cal Q}\lra {\cal
C}_{\cal Q}$ has the following property
\begin{Prop}
\label{vers} Let $\tau\colon {\cal S}:=\Spec(A)\lra \n$ be an
$\n$-scheme where $A$ is a local Artin algebra. Suppose that there
is a flat family $\sigma_{\cal S}\colon {\cal C}_{\cal S}\lra
{\cal S}$ of semistable curves over ${\cal S}$ together with an
${\cal S}$-morphism $\pi_{\cal S}\colon {\cal C}_{\cal S}\lra
\tau^* {\cal C}_\n$. Suppose that the closed point $s$ of ${\cal
S}$ maps to the origin of $\n$ and that $\pi_{\cal S|\sigma_{\cal
S}^{-1}(s)}$ equals the map $\pi$.
\par
Then, there is an $\n$-morphism $\psi\colon {\cal S}\lra {\cal
Q}$, such that ${\cal C}_{\cal S}$ is over $\tau^*{\cal C}_{\n} =
\psi^* {\cal C}_{\cal Q}$ isomorphic to $\psi^*\widetilde{\cal
C}_{\cal Q}$.
\end{Prop}
\noindent\it Proof\rm.
First note that the homomorphism
$$
H^1\bigl(\cal Hom({\Oh}mega_C^1, {\Oh}_C)\bigr)\lra H^1\bigl(\cal
Hom(\pi^*{\Oh}mega_{C_0}^1, {\Oh}_C)\bigr)
$$
is injective. In fact, as the computations in \cite{Gie} used for
proving the analogous statement (Corollary~4.4) are completely
local, they apply to our situation, too. The rest of the proof
may now be copied from \cite{Gie}, proof of Proposition~4.5. \qed
\subsubsection{Proof of Theorem~\ref{fibres}}
Let ${\frak H}^\dagger(g;\chi,r)$ be as in
Section~\ref{smoothsection}. There is a morphism
$$
\kappa^\dagger\colon\ {\frak H}^\dagger(g;\chi,r)\lra \ol{\frak
M}_g,
$$
and we define
$$
{\frak H}_{C_0}\q:=\q \kappa^{\dagger-1}\bigl([C_0]\bigr).
$$
Let $\sigma_{\frak H_{C_0}}\colon {\cal C}_{\frak H_{C_0}}\lra
{\frak H}_{C_0}$ be the restriction of the universal family. A
suitably high power of the relative dualizing sheaf $\omega_{\cal
C_{\frak H_{C_0}}/{\frak H}_{C_0}}$ will yield a morphism
$$
\begin{CD}
{\cal C}_{\frak H_{C_0}} @>\iota >> {\cal P}
\\
@V \sigma_{\frak H_{C_0}} VV @VVV
\\
{\frak H_{C_0}} @= {\frak H_{C_0}},
\end{CD}
$$
where ${\cal P}$ is some projective bundle. The image of $\iota$
is a flat family of stable curves, all of which are isomorphic to
$C_0$. As $C_0$ does not have any automorphisms, this family is
trivial. Let
$$
\phi_{\frak H_{C_0}}\colon\  {\cal C}_{{\frak H_{C_0}}} \lra
C_0\times {\frak H_{C_0}}
$$
be the induced morphism. Let $x\in {\frak H_{C_0}}$ be a point and
${\frak U}$ its formal neighborhood. Denote the fibre of the
family ${\cal C}_{{\frak H_{C_0}}}$ over $x$ by $C$, and let
$\tau\colon {\frak U}\lra \n$ be the constant map to the origin.
Finally, define
$$
\sigma_{\frak U}\colon\  {\cal C}_{\frak U} \lra {\frak U}
$$
as the restriction of the family ${\cal C}_{\frak H_{C_0}}$. By
Proposition~\ref{vers}, there is a morphism
$$
\psi\colon\  {\frak U}\lra {\cal X}\hookrightarrow {\cal Q}.
$$
Our observation in Remark~\ref{deformations}, i), implies that
the morphism $\psi$ is smooth. Therefore, ${\frak H_{C_0}}$ has
only analytic normal crossing singularities at $x$. Finally, as
$C_0$ does not have any automorphisms, the quotient morphism
$$
{\frak H_{C_0}}\cap {\frak H}^{0}(g;\chi,r)^{\rm s}\lra
\kappa_g^{-1}\bigl([C_0]\bigr)\cap {\frak H\frak C}^{\rm
s}(g;\chi,r)
$$
is a principal $\PGL(V^{\chi_l})$-bundle. This proves the theorem.
\qed
\begin{Rem}
If the automorphism group of $C_0$ is non-trivial, the same
arguments show that the fibre
$\kappa_g^{-1}\bigl([C_0]\bigr)\cap{\frak H\frak C}^{\rm
s}(g;\chi,r)$ is the quotient of a variety with analytical normal
crossings by the automorphism group of $C_0$. However, even if
$C_0$ is a smooth curve, the action of the group ${\rm Aut}(C_0)$
on the moduli space of semistable bundles has not been thoroughly
studied, so far. We refer the reader to the paper \cite{Aut} for
information concerning the action of a single automorphism.
\end{Rem}
\subsection{The moduli stacks}
\label{Stack}
In the following, let $\ul{\rm Schemes}_\C$ be the category of
schemes of finite type over $\C$, viewed as a 2-category, and
$\ul{\rm Groupoids}$ the 2-category of groupoids, that is the
2-category whose objects are groupoids, i.e., categories in which
all morphisms are isomorphisms, the 1-morphisms are functors and
the 2-morphisms are natural transformations between functors.
\par
As usual, defining presheaves of groupoids and establishing
isomorphisms between them involves many schemes characterized by
some universal property, such as fibre products. However, these
schemes will be defined only up to canonical isomorphy and there
is no equivalence relation which compensates for this. Thus, we
have to fix \sl a priori \rm a representative for every such
isomorphy class. In the following, we assume to have done this.
\par
Next, we introduce the 2-functors
$$
{\cal H\cal C}^{\rm (s)s}_{g/\chi/r}\colon\  \ul{\rm
Schemes}_\C\lra \ul{\rm Groupoids}.
$$
For any scheme $S$ of finite type, the objects of ${\cal H\cal
C}^{\rm (s)s}_{g/\chi/r}(S)$ are families $({\cal C}_{S},\E_{S})$
of H-(semi)stable vector bundles as before, the morphisms between
$({\cal C}_S^\p,\E_S^\p)$ and $({\cal C}_S,\E_S)$ are pairs
$(\phi_S,\psi_S)$, consisting of an $S$-isomorphism $\phi_S\colon
{\cal C}_S^\p\lra {\cal C}_S$ and an isomorphism $\psi_S\colon
\E^\p_S\lra \phi_S^*\E_S$. For any morphism $f\colon T\lra S$
pullback of families defines a natural transformation
$$
{\cal H\cal C}^{\rm (s)s}_{g/\chi/r}(f)\colon\  {\cal H\cal
C}^{\rm (s)s}_{g/\chi/r}(S)\lra {\cal H\cal C}^{\rm
(s)s}_{g/\chi/r}(T).
$$
On the other hand, there are the quotient stacks $\bigl[\frak
H^{\rm ss} /\GL(V^{\chi_l})\bigr]$ and $\bigl[\frak H^{\rm s}
/\GL(V^{\chi_l})\bigr]$. Here, ${\frak H}^{\rm (s)s}$ is the open
part of ${\frak H}^0(g;\chi_l,r)$ which parameterizes the
(semi)stable objects. For any scheme $S$, the objects of
$\bigl[\frak H^{\rm (s)s} /\GL(V^{\chi_l})\bigr]$ are pairs
$\bigl(\theta_S\colon {\cal P}\lra S, \eta_S\colon {\cal
P}\lra\frak H^{\rm (s)s}\bigr)$ where $\theta_S\colon {\cal
P}\lra S$ is a principal $\GL(V^{\chi_l})$-bundle and $\eta_S$ is
an equivariant morphism. One has a natural notion of isomorphism
and, as before, pullback defines the functor associated with a
morphism $f\colon T\lra S$.
\begin{Thm}
The presheaves ${\cal H\cal C}^{\rm (s)s}_{g/\chi/r}$ and
$\bigl[\frak H^{\rm (s)s} /\GL(V^{\chi_l})\bigr]$ are isomorphic.
\end{Thm}
\begin{proof} The assertions amount to prove that, for every
scheme $S$, the groupoids ${\cal H\cal C}^{\rm
(s)s}_{g/\chi/r}(S)$ and $\bigl[\frak H^{\rm (s)s}
/\GL(V^{\chi_l})\bigr](S)$ are equivalent.
\par
First, let $\bigl(\theta_S\colon {\cal P}\lra S, \eta_S\colon
{\cal P}\lra\frak H^{\rm (s)s}\bigr)$ be an object of $\bigl[\frak
H^{\rm (s)s} /\GL(V^{\chi_l})\bigr](S)$. Then, by means of
pullback, the morphism $\eta_S$ yields a
$\GL(V^{\chi_l})$-invariant, ${\cal P}$-flat family of semi\-stable
curves of genus $g$
$$
{\cal C}_{\cal P}\hookrightarrow {\cal P}\times {\frak G}
$$
and a $\GL(V^{\chi_l})$-linearized vector bundle $\E_{\cal P}$ on
${\cal C}_{\cal P}$. Now, as $S$ is the geometric quotient of
${\cal P}$ by the $\GL(V^{\chi_l})$-action, the same arguments
which were used in Section~\ref{UnivFam} show that we have the
curve $\pi_S\colon {\cal C}_S\lra S$. This time, as there are no
stabilizers present, every fibre of $\pi_S$ is indeed a
semistable curve of genus $g$. We have to check that the family
${\cal C}_S$ is indeed $S$-flat. For this, choose a Zariski-open
set $U\subset S$ over which the principal bundle ${\cal P}$ is
trivial (this is possible, since we are dealing with
$\GL(V^{\chi_l})$.) Set $V:=U\times\GL(V^{\chi_l})$. We will show
that ${\cal C}_{S|V}$ is in fact $\GL(V^{\chi_l})$-equivariantly a
product ${\cal C}_U\times \GL(V^{\chi_l})$ for some $U$-flat
family ${\cal C}_U$. This clearly settles the affair. By the
$\GL(V^{\chi_l})$-equivariance of $\theta_S$, we have the
commutative diagram
$$
\begin{array}{cccccc}
(x,g,h) &\in & U\times \GL(V^{\chi_l})\times \GL(V^{\chi_l})
&\stackrel{\eta_{S|V}\times \id}{\lra} &{\frak H}^{\rm
(s)s}\times \GL(V^{\chi_l})
\\
\downarrow &&\downarrow &&\downarrow
\\
(x, g\cdot h)&\in& U\times \GL(V^{\chi_l})
&\stackrel{\eta_{S|V}}{\lra} & {\frak H}^{\rm (s)s}&.
\end{array}
$$
Now, consider the map
\begin{eqnarray*}
V=U\times \GL(V^{\chi_l}) &\lra & \bigl(U\times
\GL(V^{\chi_l})\bigr)\times \GL(V^{\chi_l})
\\
(x,g) &\lma & (x,{\id}_{V^{\chi_l}}, g).
\end{eqnarray*}
Define $\eta^0\colon U\lra {\frak H}^{\rm (s)s}$ by $\eta^0(x):=
\eta_{S|V}(x,{\id}_{V^{\chi_l}})$. The content of the diagram
before may then be summarized by the suggestive formula
$$
\eta_{S|V}(x,g)\q=\q \eta^0(x)\cdot g.
$$
Finally, the morphism $\widetilde{\eta}\colon V \lra {\frak
H}^{\rm (s)s}$, $(x,g)\lma \eta^0(x)\cdot g$, is by definition of
the group action obtained in the following manner: Let $({\cal
C}_U,\E_U)$ be the family induced by the morphism $\eta^0$. Note
that we get even a quotient
$$
q_U\colon V^{\chi_l}\otimes{\Oh}_{\cal C_U} \lra \E_U.
$$
Let
$$
\Gamma\colon V\otimes{\Oh}_{\GL(V^{\chi_l})} \lra
V\otimes{\Oh}_{\GL(V^{\chi_l})}
$$
be the tautological automorphism. Then, define the following
quotient on ${\cal C}_U\times \GL(V^{\chi_l})$
$$
\begin{CD}
q_{V}\colon\  V\otimes{\Oh}_{\cal C_U\times \GL(V^{\chi_l})} @>
\pi_{\GL(V^{\chi_l})}^*(\Gamma) >> V\otimes{\Oh}_{\cal C_U\times
\GL(V^{\chi_l})} @> \pi_{\cal C_U}^*(q_U) >> \pi^*_{\cal C_U}\E_U.
\end{CD}
$$
This quotient defines an embedding ${\cal C}_U\times
\GL(V^{\chi_l})\lra V\times {\frak G}$, and the resulting
morphism $V\lra {\frak H}^{\rm (s)s}$ is just
$\widetilde{\eta}=\eta_{S|V}$. In particular, ${\cal C}_{S|V}$ is
$\GL(V^{\chi_l})$-equivariantly isomorphic to $ {\cal C}_U\times
\GL(V^{\chi_l})$, as asserted.
By Kempf's descent lemma, the bundle $\E_{\cal P}$ descends to
${\cal C}_S$, so that $({\cal C}_S,\E_S)$ is an object of ${\cal
H\cal C}^{\rm (s)s}_{g/\chi/r}(S)$. An isomorphism in $\bigl[\frak
H^{\rm (s)s} /\GL(V^{\chi_l})\bigr](S)$ will clearly lead to a
unique isomorphism in ${\cal H\cal C}^{\rm(s)s}_{g/\chi/r}(S)$.
\par
Now, suppose we are given a scheme $S$ and a family $(\pi_S\colon
{\cal C}_S\lra S,\E_S)$ of H-(semi)stable vector bundles. Then,
we know that $\theta_S\colon \cal Isom(V^{\chi_l}\otimes{\Oh}_S,
\pi_{S*}\E_S)\lra S$ is a principal $\GL(V^{\chi_l})$-bundle. On
${\cal P}:=\cal Isom(V^{\chi_l}\otimes{\Oh}_S, \pi_{S*}\E_S)$, there
is the tautological isomorphism
$$
\tau_{\cal P}\colon\  V^{\chi_l}\otimes{\Oh}_{\cal P}\lra
\theta_S^*\pi_{S*}\E_S.
$$
Now, form the cartesian diagram
$$
\begin{CD}
{\cal C}_{\cal P} @> \psi>> {\cal C}_S
\\
@V\pi_{\cal P} VV @VV \pi_S V
\\
\cal P @>\theta_S >> S.
\end{CD}
$$
By flat base change
$$
\theta_S^*\pi_{S*}\E_S\q\cong\q \pi_{\cal P*}\psi^*\E_S.
$$
If we set $\E_{\cal P}:=\psi^* \E_S$, then
$$
V^{\chi_l}\otimes{\Oh}_{\cal C_{\cal P}}\stackrel{\pi_{\cal
P}^*\tau_{\cal P}}{\lra} \pi_{\cal P}^*\pi_{\cal P*}\E_{\cal P}
\stackrel{\rm ev}{\lra} \E_{\cal P}
$$
defines a morphism $\cal P\lra {\frak H}^{\rm (s)s}$ which is by
construction $\GL(V^{\chi_l})$-equivariant. Again, isomorphisms
in the category ${\cal H\cal C}^{\rm(s)s}_{g/\chi/r}(S)$ will lead
canonically to isomorphisms in the groupoid $\bigl[\frak H^{\rm
(s)s} /\GL(V^{\chi_l})\bigr](S)$.
\par
The two operations just introduced clearly establish the desired
equivalence of categories. 
\end{proof}
\par
Now, by the results of Section~\ref{smoothsection}, we know that
the schemes ${\frak H}^{\rm ss}$ and ${\frak H}^{s}$ are smooth.
Moreover, the quotient map ${\frak H}^{\rm (s)s}\lra \bigl[\frak
H^{\rm (s)s}/\GL(V^{\chi_l})\bigr]$ is smooth, whence $\bigl[\frak
H^{\rm ss}/\GL(V^{\chi_l})\bigr]$ is a smooth Artin stack and
$\bigl[{\frak H}^{\rm s}/\GL(V^{\chi_l})\bigr]$ is a smooth
Deligne-Mumford stack.
\begin{Cor}
The Hilbert compactification ${\cal H\cal C}^{\rm ss}_{g/\chi/r}$
is a smooth Artin stack, and its open substack ${\cal H\cal
C}^{\rm s}_{g/\chi/r}$ is a smooth Deligne-Mumford stack.
\end{Cor}

\end{document}